\documentclass[reqno,10pt]{amsart}
\usepackage{xcolor}
\usepackage{amssymb}
\usepackage[all,cmtip]{xy}
\newtheorem{MainTheorem}{Theorem}
\newtheorem{Proposition}{Proposition}[section]
\newtheorem{Definition}[Proposition]{Definition}
\newtheorem{Lemma}[Proposition]{Lemma}

\newtheorem{Corollary}[Proposition]{Corollary}
\newtheorem{Remark}[Proposition]{Remark}

\DeclareMathOperator{\Val}{Val}

\DeclareMathOperator{\PD}{PD}
\DeclareMathOperator{\vol}{vol}
\DeclareMathOperator{\conv}{conv}
\DeclareMathOperator{\spt}{spt}
\DeclareMathOperator{\Dens}{Dens}
\DeclareMathOperator{\Span}{Span}
\DeclareMathOperator{\dist}{dist}
\DeclareMathOperator{\Image}{Im}
\DeclareMathOperator{\WF}{WF}
\DeclareMathOperator{\ori}{or}
\DeclareMathOperator{\singsupp}{singsupp}
\DeclareMathOperator{\GL}{GL}
\newcommand{\R}{\mathbb{R}}

\title[Generalized valuations and the polytope algebra ]{Generalized translation invariant valuations and the polytope
algebra}
\author{Andreas Bernig}
\author{Dmitry Faifman}
 
\email{bernig@math.uni-frankfurt.de}
\email{dfaifman@ihes.fr}

\address{Institut f\"ur Mathematik, Goethe-Universit\"at Frankfurt,
Robert-Mayer-Str. 10, 60054 Frankfurt, Germany}
\address{Institut des Hautes \'Etudes Scientifiques, 35 Route de Chartres, 91440 Bures-sur-Yvette, France}

\thanks{A. B. was supported by DFG grants BE 2484/3-1 and BE 2484/5-1. D. F. was partially supported by ISF grant
1447/12.\\ AMS 2010 {\it Mathematics subject
classification}: 
52B45; 
53C65
}

\begin{document}

\begin{abstract}
We study the space of generalized translation invariant valuations on a finite-dimensional vector space and
construct
a partial convolution which extends the convolution of smooth translation
invariant valuations. Our main theorem is that McMullen's polytope algebra is a subalgebra of the (partial)
convolution algebra of generalized translation invariant valuations. More precisely, we show that the 
polytope algebra embeds injectively into the space of generalized
translation invariant valuations and that for polytopes in general position, the convolution is defined
and corresponds to the product in the polytope algebra.  
\end{abstract}

\maketitle 

\section{Introduction}

Let $V$ be an $n$-dimensional vector space, $V^*$ the dual vector space, $\mathcal{K}(V)$ the set of 
non-empty compact convex
subsets in $V$, endowed with the topology induced by the Hausdorff metric for an arbitrary Euclidean structure on $V$,
and $\mathcal{P}(V)$ the set of polytopes in $V$. A valuation is a
map $\mu:\mathcal{K}(V) \to \mathbb{C}$ such that 
\begin{displaymath}
 \mu(K \cup L)+\mu(K \cap L)=\mu(K)+\mu(L)
\end{displaymath}
whenever $K,L, K \cup L \in \mathcal{K}(V)$. Continuity of valuations will be with respect to the Hausdorff topology.

Examples of valuations are measures, the intrinsic volumes (in particular the Euler characteristic $\chi$) and mixed
volumes. 

Let $\Val(V)$ denote the (Banach-)space of continuous, translation invariant valuations. It was the object of
intensive research during the last few years, compare \cite{alesker_mcullenconj01, alesker_fourier,
alesker_bernig_schuster, alesker_faifman, bernig_aig10, bernig_broecker07, bernig_fu06, bernig_fu_hig, bernig_hug,
fu_barcelona} and
the references therein. 

Valuations with values in semi-groups other than $\mathbb{C}$ have also attracted a lot of interest. We only
mention the recent papers \cite{abardia12, abardia_bernig, bernig_fu_solanes, haberl10,
hug_schneider_localtensor, schneider13, schuster10, schuster_wannerer, wannerer_area_measures, wannerer_unitary_module} 
to give a flavor on this active research area.

Of particular importance is the class of the so-called smooth valuations. The importance of this class stems from the fact that it admits various algebraic structures, which include two bilinear pairings, known as product and convolution, and a Fourier-type duality interchanging them. These algebraic structures are closely related to important notions from convex and integral geometry, such as the Minkowski sum, mixed volumes, and kinematic formulas. This emerging new theory is known as algebraic integral geometry \cite{bernig_aig10, fu_barcelona}.

A different, more classical type of algebraic object playing an important role in convex geometry is McMullen's algebra of polytopes. 
In this paper, we show how McMullen's algebra fits into the framework of algebraic integral geometry. More precisely, we show that McMullen's algebra can be embedded as a subalgebra of the space of generalized valuations, which is, roughly speaking, the dual space of smooth valuations.

Let us now give the necessary background required to state our main theorems. 

The group $\GL(V)$ acts in the natural way on $\Val(V)$. The dense subspace of $\GL(V)$-smooth vectors in $\Val(V)$ is
denoted by $\Val^\infty(V)$. It carries a Fr\'echet topology which is finer than the induced topology.  

In \cite{bernig_fu06}, a convolution product on $\Val^\infty(V) \otimes \Dens(V^*)$ was constructed. Here
 and in the following, $\Dens(W)$
denotes the $1$-dimensional space of densities on a linear space $W$. Note that $\Dens(V) \otimes \Dens(V^*) \cong
\mathbb{C}$: if $\vol$ is
any choice of Lebesgue measure on $V$, and $\vol^*$ the corresponding dual measure on $V^*$, then $\vol \otimes \vol^*
\in \Dens(V) \otimes \Dens(V^*)$ is independent of the choice of $\vol$. If
$\phi_i(K)=\vol(K+A_i) \otimes \vol^*$ with smooth compact strictly convex bodies $A_1,A_2$, then $\phi_1 *
\phi_2(K)=\vol(K+A_1+A_2) \otimes \vol^*$. By Alesker's  proof \cite{alesker_mcullenconj01} of
McMullen's conjecture,
linear
combinations of such valuations are dense in the space of all smooth valuations. The convolution
extends by bilinearity and continuity to $\Val^\infty(V) \otimes \Dens(V^*)$.

By \cite{bernig_fu06}, the convolution product is closely related to additive
kinematic formulas. It was recently used in the study of unitary kinematic formulas \cite{bernig_fu_hig}, local unitary
additive kinematic formulas \cite{wannerer_area_measures, wannerer_unitary_module} and kinematic formulas for tensor
valuations \cite{bernig_hug}.
 
 In this paper, we will extend the convolution to a (partially defined) convolution on the space of generalized
translation invariant valuations.

\begin{Definition}
Elements of the space 
 \begin{displaymath}
  \Val^{-\infty}(V):=\Val^\infty(V)^* \otimes \Dens(V)
 \end{displaymath}
are called generalized translation invariant valuations. 
\end{Definition}

By the Alesker-Poincar\'e duality \cite{alesker04_product}, $\Val^\infty(V)$ embeds in $\Val^{-\infty}(V)$ as a dense
subspace. More generally, it follows from \cite[Proposition 8.1.2]{alesker_fourier} that $\Val(V)$ embeds in
$\Val^{-\infty}(V)$, hence we have the inclusions
\begin{displaymath}
\Val^\infty(V) \subset \Val(V) \subset \Val^{-\infty}(V).
\end{displaymath}

Generalized translation invariant valuations were introduced and studied in the recent paper \cite{alesker_faifman}.
Note that another notion of generalized valuation was introduced by Alesker in \cite{alesker_val_man1,
alesker_val_man2, alesker_val_man4, alesker_val_man3}. In the next section, we will construct a natural
isomorphism between the space of translation invariant generalized valuations in Alesker's sense and the space
of generalized translation invariant valuations in the sense of the above definition.

Given a polytope $P$ in $V$, there is an element $M(P) \in \Val^{-\infty}(V) \otimes
\Dens(V^*) \cong
\Val^\infty(V)^*$ defined by 
\begin{displaymath}
 \langle M(P),\phi\rangle=\phi(P). 
\end{displaymath}

Let $\Pi(V)$ be McMullen's polytope algebra \cite{mcmullen_polytope_algebra}. As a vector space, $\Pi(V)$ is
generated by all symbols
$[P]$, where $P$ is a polytope in $V$, modulo the relations $[P] \equiv [P+v], v \in V$, and $[P \cup Q]+[P \cap Q]=[P]
+ [Q]$ whenever $P,Q, P \cup Q$ are polytopes in $V$. The product is defined by $[P] \cdot [Q]:=[P+Q]$. 

The map $M:\mathcal{P}(V) \to \Val^{-\infty}(V) \otimes \Dens(V^*)$ extends to a linear map 
\begin{displaymath}
 M: \Pi(V) \to \Val^{-\infty}(V) \otimes \Dens(V^*).
\end{displaymath}

Our first main theorem shows that McMullen's polytope algebra is a subset of
$\Val^{-\infty}(V)\otimes
\Dens(V^*)$. 

\begin{MainTheorem} \label{mainthm_injection_mcmullen}
The map $M:\Pi(V)\to \Val^{-\infty}(V)\otimes \Dens(V^*)$ is injective. 
Equivalently, the elements of
$\Val^\infty(V)$ separate the elements of $\Pi(V)$.
\end{MainTheorem}

In Section \ref{sec_partial_conv} we will introduce a notion of transversality of generalized translation invariant
valuations. Our second main theorem is the following. 

\begin{MainTheorem} \label{mainthm_convolution}
There exists a partial convolution product $*$ on $\Val^{-\infty}(V) \otimes \Dens(V^*)$ with the following properties:
\begin{enumerate}
 \item If $\phi_1,\phi_2 \in \Val^{-\infty}(V) \otimes \Dens(V^*)$ are transversal, then $\phi_1 * \phi_2 \in
\Val^{-\infty}(V) \otimes \Dens(V^*)$ is defined. 
\item If $\phi_1 * \phi_2$ is defined and $g \in \GL(V)$, then $(g_* \phi_1) * (g_* \phi_2)$ is defined
and equals $g_*(\phi_1 * \phi_2)$.
\item Whenever the convolution is defined, it is bilinear, commutative, associative and of degree $-n$. 
\item The restriction to the subspace $\Val^\infty(V) \otimes \Dens(V^*)$ is the convolution product from
\cite{bernig_fu06}. 
\item If $x,y$ are elements in $\Pi(V)$ in general position, then $M(x),M(y)$ are transversal in $\Val^{-\infty}(V)
\otimes \Dens(V^*)$ and 
\begin{displaymath}
 M(x \cdot y)=M(x) * M(y). 
\end{displaymath}
\end{enumerate}
\end{MainTheorem}

Stated otherwise, the maps in the diagram 
\begin{displaymath}
  \xymatrix{\Val^\infty(V) \otimes \Dens(V^*) \ar@{^{(}->}[r] & \Val^{-\infty}(V) \otimes \Dens(V^*) & \Pi(V)
\ar@{_{(}->}[l]}
\end{displaymath}
have dense images and are compatible with the (partial) algebra structures. 

{\bf Remarks:} 
\begin{enumerate}
\item In Prop. \ref{prop_continuity} we will show that, under some technical conditions in terms of wave fronts,
the convolution on generalized translation invariant valuations from Theorem \ref{mainthm_convolution} is the unique
jointly sequentially continuous extension of the convolution product on smooth translation invariant valuations.
\item In \cite{alesker_bernig}, it was shown that the space $\mathcal{V}^{-\infty}(X)$  of generalized valuations
on a smooth manifold $X$ admits a partial product structure extending the Alesker product of smooth valuations
on $X$. If $X$ is
real-analytic, then the space $\mathcal{F}_\mathbb{C}(X)$ of $\mathbb{C}$-valued constructible functions on $X$ embeds
densely into $\mathcal{V}^{-\infty}(X)$. It was conjectured that whenever two constructible functions meet
transversally, then the
product in the sense of generalized valuations exists and equals the generalized valuation corresponding to the product
of the two functions. The relevant diagram in this case is 
\begin{displaymath}
  \xymatrix{\Val^\infty(X) \ar@{^{(}->}[r] & \Val^{-\infty}(X) & \mathcal{F}_\mathbb{C}(X),
\ar@{_{(}->}[l]}
\end{displaymath}
where both maps are injections with dense images and are (conjecturally) compatible with the partial product structure.
Theorem 5 in \cite{alesker_bernig} gives strong support for this conjecture. 
\item  The Alesker-Fourier transform from \cite{alesker_fourier} extends to an isomorphism
$\mathbb{F}:\Val^{-\infty}(V^*) \to \Val^{-\infty}(V) \otimes \Dens(V^*)$, compare \cite{alesker_faifman}. Another
natural partially defined convolution on $\Val^{-\infty}(V) \otimes \Dens(V^*)$ would be 
\begin{displaymath}
\psi_1 * \psi_2:=\mathbb{F} \left(\mathbb{F}^{-1}\psi_1 \cdot \mathbb{F}^{-1} \psi_2\right), 
\end{displaymath}
where the dot is the partially defined product on $\Val^{-\infty}(V^*)$ from \cite{alesker_bernig}. It seems natural to expect that this convolution coincides with the one from Theorem \ref{mainthm_convolution}, but we do not have a proof of this fact. 
\end{enumerate}

\subsection*{Plan of the paper}

In the next section, we introduce and study the space of generalized translation invariant valuations and
explore its relation to generalized valuations from Alesker's theory. In Section \ref{sec_embedding_polytopes} we show
that McMullen's polytope algebra embeds into the space of generalized translation invariant valuations.  A partial
convolution structure on this space is constructed in Section
\ref{sec_partial_conv}. In Section \ref{sec:Appendix} we construct a
certain current on the sphere which is related to the volumes of spherical joins, and can be viewed as a generalization of the Gauss area formula for the plane. This current also plays a major role in
the proof of the second main theorem. Its construction is based on geometric measure theory and is independent of
the rest of the paper. Finally, Section \ref{sec_compatibility} is devoted to the proof of the fact that the embedding of the
polytope algebra is compatible with the two convolution (product) structures. 
\subsection*{Acknowledgements}
 We wish to thank Semyon Alesker for multiple fruitful discussions and Thomas Wannerer for useful remarks on a first
draft of this paper.

\section{Preliminaries}
\label{sec_prels}
 
In this section $X$ will be an oriented $n$-dimensional smooth manifold and $S^*X$ the cosphere bundle over $X$. It consists of
all pairs
$(x,[\xi])$ where $x \in X, \xi \in T^*_xX, \xi \neq 0$ and where the equivalence relation is defined by $[\xi]=[\tau]$
if and only if $\xi=\lambda \tau$ for some $\lambda>0$.

The projection onto $X$ is denoted by $\pi:S^*X \to X,
(x,\xi) \mapsto x$. The antipodal map $s:S^*X \to S^*X$ is defined by $(x,[\xi]) \mapsto (x,[-\xi])$.

The push-forward map (also called fiber integration) $\pi_*:\Omega^k(S^*X) \to \Omega^{k-(n-1)}(X)$ satisfies 
\begin{displaymath}
\int_{S^*X} \pi^*\gamma \wedge \omega=\int_X \gamma \wedge \pi_*\omega, \quad \gamma \in \Omega_c^{2n-k-1}(X). 
\end{displaymath}

If $V$ is a vector space, then $S^*V \cong V \times \mathbb P_+(V^*)$, where $\mathbb P_+(V^*):=V^*
\setminus \{0\}/ \R_+$ is the sphere in $V^*$. Moreover, if $V$ is 
a Euclidean vector space of dimension $n$, we identify $S^*V$ and $SV=V \times S^{n-1}$ and write $\pi=\pi_1:SV \to V,
\pi_2:SV \to S^{n-1}$ for the two projections. 

\subsection{Currents}

Let us recall some terminology from geometric measure theory. We refer to \cite{federer_book, morgan_book} for more
information. 

The space of $k$-forms on $X$ is denoted by $\Omega^k(X)$, the space of compactly
supported $k$-forms is denoted by $\Omega^k_c(X)$. Elements of the  dual space $\mathcal{D}_k(X):=\Omega^k_c(X)^*$ are
called {\it $k$-currents}. A $0$-current is also called {\it distribution}.

The boundary of a $k$-current $T \in \mathcal{D}_k(X)$ is defined by $\langle \partial
T,\phi\rangle=\langle T,d\phi\rangle, \phi \in \Omega^{k-1}_c(X)$. If $\partial T=0$, $T$ is called a {\it cycle}. 
If $T \in \mathcal{D}_k(X)$ and $\omega \in \Omega^l(X), l \leq k$, then the current $T \llcorner \omega \in
\mathcal{D}_{k-l}(X)$ is defined by $\langle T \llcorner \omega,\phi\rangle:=\langle T, \omega \wedge \phi\rangle$.

If $f:X \to Y$ is a smooth map between smooth manifolds $X,Y$ and $T \in \mathcal{D}_k(X)$ such that $f|_{\spt
T}$ is proper, then the push-forward $f_*T \in \mathcal{D}_k(Y)$ is defined by $\langle f_*T,\phi\rangle:=\langle
T,\zeta f^*\phi\rangle$, where $\zeta \in C^\infty_c(X)$ is equal to $1$ in a neighborhood of $\spt T \cap \spt
f^*\omega$. It is easily checked that 

\begin{equation}\label{eq_boundary_vs_differential}
 \partial ([[X]] \llcorner \omega)=(-1)^{\deg\omega+1} [[X]] \llcorner (d\omega)
\end{equation}
and
\begin{equation} \label{eq_push_forward_forms_currents}
 \pi_*\left([[S^*X]] \llcorner \omega\right)=(-1)^{(n+1)(\deg\omega+1)} [[X]] \llcorner \pi_*\omega.
\end{equation}

Every oriented submanifold $Y \subset X$ of dimension $k$ induces a $k$-current $[[Y]]$ such that $\langle
[[Y]],\phi\rangle=\int_Y \phi$. By Stokes' theorem, $\partial [[Y]]=[[\partial Y]]$. A {\it smooth current} is a current
of the form $[[X]] \llcorner \omega \in \mathcal{D}_{n-k}(X)$ with  $\omega \in \Omega^k(X)$.

If $X$ and $Y$ are smooth manifolds, $T \in \mathcal{D}_k(X), S \in \mathcal{D}_l(Y)$, then there is a unique
current $T \times S \in \mathcal{D}_{k+l}(X \times Y)$ such that $\langle T \times S,
\pi_1^*\omega \wedge \pi_2^*\phi\rangle=\langle T,\omega\rangle \cdot \langle S,\phi\rangle$, for all $\omega \in
\Omega^k(X), \phi \in \Omega^l(Y)$. Here $\pi_1,\pi_2$ are the projections from $X \times Y$ to $X$ and $Y$
respectively. 

If $T=[[X]] \llcorner \omega, S=[[Y]] \llcorner \phi$, then 
\begin{equation} \label{eq_product_smooth_currents}
 T \times S=(-1)^{(\dim X-\deg \omega)\deg \phi}[[X \times Y]] \llcorner (\omega \wedge \phi).
\end{equation}

The boundary of the product is given by 
\begin{equation} \label{eq_boundary_product}
 \partial(T \times S)=\partial T \times S+(-1)^k T \times \partial S,
\end{equation}
compare \cite[4.1.8]{federer_book}.

If $X$ is a Riemannian manifold, the {\it mass} of a current $T \in \mathcal{D}_k(X)$ is 
\begin{displaymath}
 \mathbf{M}(T):=\sup \{\langle T,\phi\rangle: \phi \in \Omega^k_c(X), \|\phi(x)\|^* \leq 1, \forall x \in X\},
\end{displaymath}
where $\|\cdot\|^*$ denotes the comass norm. 

Currents of finite mass having a boundary of finite mass are called {\it normal currents}. 

The flat norm of $T$ is defined by 
\begin{displaymath}
 \mathbf{F}(T):=\sup\{\langle T,\phi\rangle: \phi \in \Omega^k_c(X), \|\phi(x)\|^* \leq 1, \|d\phi(x)\|^* \leq 1,
\forall x \in X\}. 
\end{displaymath}
If $X$ is compact, then the $\mathbf{F}$-closure of the space of normal $k$-currents is the space of {\it real flat
chains}.

\subsection{Wave fronts}

We refer to \cite{guillemin_sternberg77} and \cite{hoermander_pde1} for the general theory of wave fronts and
its applications. For the reader's convenience and later reference, we will recall some basic definitions and
some fundamental properties of wave fronts, following \cite{hoermander_pde1}.

First let $X$ be a linear space of dimension $n$, $T$ a distribution on $X$.

The cone $\Sigma(T) \subset X^*$ is defined as the closure of the complement of the set of all $\eta \in
X^*$ such that for all $\xi$ in a conic neighborhood of $\eta$ we have 
\begin{displaymath}
 \|\hat T(\xi)\| \leq C_N(1+\|\xi\|)^{-N}, \quad N \in \mathbb{N}
\end{displaymath}
(with constants $C_N$ only dependent on $N$ and the chosen neighborhood). Here $\hat T$ denotes the Fourier
transform
of $T$, and the norm is taken with respect to an arbitrary scalar product on $X^*$. 

Next, for an affine space $X$ and a point $x \in X$, the set $\Sigma_x(T) \subset T_x^* X$ is defined by 
\begin{displaymath}
 \Sigma_x(T):=\bigcap_\phi \Sigma(\phi T),
\end{displaymath}
where $\phi$ ranges over all compactly supported smooth functions on $X$ with $\phi(x) \neq 0$. Note that one uses the
canonic identification $T^*_xX=X^*$.

The {\it wave front} set of $T$ is by definition 
\begin{displaymath}
 \WF(T):=\{(x,[\xi]) \in S^*X: \xi \in \Sigma_x(T)\}.
\end{displaymath}

The set $\singsupp(T):=\pi(\WF(T)) \subset X$ is called {\it singular support}. 

Let $(x_1,\ldots,x_n)$ be coordinates on $X$. Given a current $T \in \mathcal{D}_k(X)$, we may write
\begin{displaymath}
T=\sum_{\substack{I=(i_1,\ldots,i_k)\\1 \leq i_1<i_2<\cdots<i_k \leq n}}
T_I \frac{\partial}{\partial x_{i_1}} \wedge \ldots \wedge \frac{\partial}{\partial x_{i_k}}
\end{displaymath}
with distributions $T_I$. Then the wave front of $T$ is defined as $\WF(T):=\bigcup_I \WF(T_I)$. 

A current $T$ is smooth, i.e. given by integration against a smooth differential form,
if and only if $\WF(T)=\emptyset$.

\begin{Definition} \label{def_dkgamma}
Let $\Gamma \subset S^*X$ be a closed set. Then  we set
\begin{displaymath}
 \mathcal{D}_{k,\Gamma}(X):=\{T \in \mathcal{D}_k(X): \WF(T) \subset \Gamma\}.
\end{displaymath}
A sequence $T_j \in  \mathcal{D}_{k,\Gamma}(X)$ converges to $T \in  \mathcal{D}_{k,\Gamma}(X)$ if (writing $T_j, T$ as above) $T_{j,I} \to T_I$ weakly in the sense of distributions and for each compactly supported function $\phi \in C^\infty_c(X)$ and each closed cone $A$ in $X^*$ such that $\Gamma \cap (\spt \phi \times A)=\emptyset$ we have 
\begin{displaymath}
 \sup_{\xi \in A} |\xi|^N \left|\widehat{\phi T_{I}}(\xi)-\widehat{\phi T_{j,I}}(\xi)\right| \to 0, \quad j \to \infty
\end{displaymath}
for all $N \in \mathbb{N}$.
\end{Definition}

\begin{Proposition}[{\cite[Thm. 8.2.3]{hoermander_pde1}}] \label{prop_approx_smooth}
 Let $T \in \mathcal{D}_{k,\Gamma}(X)$. Then there exists a sequence of compactly supported smooth $k$-forms
$\omega_i \in \Omega^{n-k}(X)$ such that $[[X]] \llcorner \omega_i \to T$ in $\mathcal{D}_{k,\Gamma}(X)$. In other
words, smooth forms are dense in $\mathcal{D}_{k,\Gamma}(X)$. 
\end{Proposition}

\begin{Proposition}[{\cite[Thm. 8.2.10]{hoermander_pde1}}] \label{prop_intersection_current}
 Let $T_1 \in \mathcal{D}_{k_1}(X), T_2 \in \mathcal{D}_{k_2}(X)$ such that the following transversality condition is
satisfied: 
\begin{displaymath}
 \WF(T_1) \cap s \WF(T_2) = \emptyset.
\end{displaymath}
Then the intersection current $T_1 \cap T_2 \in \mathcal{D}_{k_1+k_2-n}(X)$ is well-defined. More precisely, if
$[[X]] \llcorner \omega_i^j \to T_j$ in $\mathcal{D}_{k_j,\WF(T_j)}(X)$ with $\omega_i^j \in \Omega^{n-k_j}(X)$,
$j=1,2$, then $[[X]] \llcorner (\omega_i^1 \wedge \omega_i^2) \to T_1 \cap T_2$ in
$\mathcal{D}_{k_1+k_2-n,\Gamma}(X)$,   where 
\begin{displaymath}
 \Gamma:=\WF(T_1) \cup \WF(T_2) \cup \left\{(x,[\xi_1+\xi_2]):(x,[\xi_1]) \in \WF(T_1), (x,[\xi_2]) \in
\WF(T_2)\right\}.
\end{displaymath}
The boundary of the intersection is given by 
\begin{equation} \label{eq_boundary_intersection}
 \partial (T_1 \cap T_2)=(-1)^{n-k_2} \partial T_1 \cap T_2+T_1 \cap \partial T_2.
\end{equation}
\end{Proposition}

The wave front of a distribution is defined locally and behaves well under coordinate changes. Using local coordinates,
one can define the wave front set $\WF(T) \subset S^*X$ for a distribution $T$ on a smooth manifold $X$. Definition \ref{def_dkgamma} and Propositions \ref{prop_approx_smooth} and \ref{prop_intersection_current} remain valid in this greater generality. 

We will need a special case of these constructions. 

\begin{Proposition} \label{prop_wf_submfld}
 Let $Y \subset X$ be a compact oriented $k$-dimensional submanifold. Then 
\begin{displaymath}
 \WF([[Y]])=N_X(Y)=\{(x,[\xi]) \in S^*X|_Y:\xi|_{T_xY}=0\}.
\end{displaymath}
\end{Proposition}

An example for Proposition \ref{prop_intersection_current} is when $T_i=[[Y_i]]$ with oriented submanifolds $Y_1,Y_2 \subset X$ intersecting
transversally (in the usual sense). Then Proposition \ref{prop_wf_submfld} implies that the transversality condition in Proposition \ref{prop_intersection_current} is satisfied, and $T_1 \cap T_2=[[Y_1 \cap Y_2]]$.

It is easily checked using \eqref{eq_product_smooth_currents}, that for currents $A_1,A_2$ on an
$n$-dimensional manifold $X$ and currents $B_1,B_2$ on an $m$-dimensional manifold $Y$, we have 
\begin{equation} \label{eq_intersection_product}
(A_1 \times B_1) \cap (A_2 \times B_2) = (-1)^{(n-\deg A_1)(m-\deg B_2)} (A_1 \cap A_2) \times (B_1 \cap B_2)
\end{equation} 
whenever both sides are well-defined.

Given a differential operator, we have $\WF(PT) \subset \WF(T)$ with equality
in case $P$ is elliptic \cite[(8.1.11) and Corollary 8.3.2]{hoermander_pde1}. In
particular, it follows that for a current on a manifold $X$, we have 
\begin{equation} \label{eq_wavefront_boundary}
 \WF(\partial T) \subset \WF(T).
\end{equation}
If $T$ is a current on the sphere $S^{n-1}$ and $\Delta$ the Laplace-Beltrami operator, then 
\begin{equation} \label{eq_wavefront_laplacian}
 \WF(\Delta T)=\WF(T). 
\end{equation}

\subsection{Valuations}

Let us now briefly recall some notions from Alesker's theory of valuations on manifolds, referring to
\cite{alesker_val_man1, alesker_val_man2, alesker_survey07, alesker_val_man4, alesker_val_man3} for more details. 
  
Let $X$ be a smooth manifold of dimension $n$, which for simplicity we suppose to be oriented. Let
$\mathcal{P}(X)$ be the space of all compact differentiable polyhedra on $X$. Given $P \in \mathcal{P}(X)$, the
conormal cycle
$N(P)$ is a Legendrian cycle in the cosphere bundle $S^*X$ (i.e. $\partial N(P)=0$ and $N(P) \llcorner \alpha=0$
where $\alpha$ is the contact form on $S^*X$). A map of the form 
\begin{displaymath}
 P \mapsto \int_{N(P)} \omega+ \int_P \gamma, \quad P \in \mathcal{P}(X), \omega \in \Omega^{n-1}(S^*X), \gamma \in
\Omega^n(X)
\end{displaymath}
is called a {\it smooth valuation} on $X$. The space of smooth valuations on $X$ is denoted by $\mathcal{V}^\infty(X)$.
It carries a natural Fr\'echet space topology. The subspace of compactly supported smooth valuations is denoted by
$\mathcal{V}_c^\infty(X)$. The valuation defined by the above equation will be denoted by $\nu(\omega,\gamma)$. 

We remark that, without using an orientation, we can still define $\nu(\omega,\phi)$, where $\omega \in
\Omega^{n-1}(S^*X) \otimes \ori(X), \phi \in \Omega^n(X) \otimes \ori(X)$. Here $\ori(X)$ is the orientation bundle over
$X$.

Elements of the space 
\begin{displaymath}
 \mathcal{V}^{-\infty}(X):=(\mathcal{V}_c^\infty(X))^* 
\end{displaymath}
are called {\it generalized valuations} \cite{alesker_val_man4}. Each compact differentiable polyhedron $P$ defines a
generalized valuation $\Gamma(P)$ by 
\begin{displaymath}
\langle \Gamma(P),\phi\rangle:=\phi(P), \quad \phi \in \mathcal{V}^\infty_c(X). 
\end{displaymath}
A smooth
valuation can be considered as a generalized valuation by Alesker-Poincar\'e duality  \cite[Thm.
6.1.1.] {alesker_val_man4}. 

We thus have injections
\begin{displaymath}
 \xymatrix{ \mathcal{V}^\infty(X) \ar@{^{(}->}[r] & \mathcal{V}^{-\infty}(X) & \mathcal{P}(X) \ar@{_{(}->}[l]}.
\end{displaymath}

By the results in \cite{bernig_broecker07} and \cite{alesker_bernig}, a generalized valuation $\phi \in
\mathcal{V}^{-\infty}(X)$ is uniquely described by a pair of currents
$E(\phi)=(T(\phi),C(\phi)) \in \mathcal{D}_{n-1}(S^*X) \times  \mathcal{D}_n(X)$ such that
\begin{equation} \label{eq_conditions_t_c}
 \partial T=0, \pi_*T=\partial C, T \text{ is Legendrian, i.e.} T \llcorner \alpha=0.
\end{equation}

Note that, in contrast to different uses of the word {\it Legendrian} in the literature, $T$ is not assumed to be
rectifiable. 

Given $(T,C)$ satisfying these conditions, we denote by $E^{-1}(T,C)$ the corresponding generalized valuation. 

If $\mu$ is a compactly supported smooth valuation on $X$, then we may represent
$\mu=\nu(\omega,\gamma)$ with compactly supported forms $\omega,\gamma$. If $E(\phi)=(T,C)$, then 
\begin{displaymath}
 \langle \phi,\mu\rangle=T(\omega)+C(\gamma).
\end{displaymath}

In particular, the generalized valuation $\Gamma(P)$ corresponding to $P \in \mathcal{P}(X)$ satisfies
\begin{equation} \label{eq_currents_for_polytope}
E(\Gamma(P))=(N(P),[[P]]).  
\end{equation}

If $\phi=\nu(\omega,\gamma)$ is smooth, then 
\begin{align}
 T(\phi) & = [[S^*X]] \llcorner s^*(D\omega+\pi^* \gamma),\label{eq_pair_currents_smooth_a}\\ 
 C(\phi) & = [[X]] \llcorner \pi_*\omega, \label{eq_pair_currents_smooth}
\end{align}
where $D$ is the Rumin operator and $s$ is the involution on $S^*X$ given by $[\xi] \mapsto [-\xi]$
\cite{alesker_bernig, rumin94}.

Let us specialize to the case where $X=V$ is a finite-dimensional vector space.  We denote by
$\mathcal{V}^\infty(V)^{tr}$ and $\mathcal{V}^{-\infty}(V)^{tr}$ the spaces of translation invariant elements.  
Alesker has shown in  \cite{alesker_val_man2} that 
\begin{displaymath}
 \mathcal{V}^\infty(V)^{tr} \cong \Val^\infty(V).
\end{displaymath}
A similar statement for generalized valuations is shown in the next proposition. 

\begin{Proposition} \label{prop_identification_vals}
 The transpose of the map 
\begin{align*}
 F: \mathcal{V}_c^\infty(V) & \to \Val^\infty(V) \otimes \Dens(V^*) \\
\mu & \mapsto \int_V \mu(\bullet + x)  d\vol(x) \otimes \vol^*
\end{align*}
induces an isomorphism 
\begin{displaymath}
 F^*: \Val^{-\infty}(V) \stackrel{\cong}{\longrightarrow} \mathcal{V}^{-\infty}(V)^{tr}.
\end{displaymath}
The diagram 
\begin{displaymath}
 \xymatrix{\Val^{-\infty}(V) \ar_\cong^{F^*}[r] & \mathcal{V}^{-\infty}(V)^{tr}\\
\Val^\infty(V) \ar^{\cong}[r] \ar@{^{(}->}[u] & \mathcal{V}^\infty(V)^{tr} \ar@{^{(}->}[u]}
\end{displaymath}
commutes and the vertical maps have dense images. 
\end{Proposition}

\proof
First we show that the diagram is commutative, i.e. that the restriction of $F^*$ to $\Val^\infty(V)$
is the
identity. 

Let $\phi\in \Val^\infty(V)$ be
a smooth valuation. We will show that for any $\mu \in \mathcal V_c^\infty(V)$ one has 
\begin{displaymath}
\langle F^*\phi,\mu\rangle_{\mathcal V^\infty(V)}=\langle \phi,\mu\rangle_{\mathcal V^\infty(V)},
\end{displaymath}
or equivalently that 
\begin{displaymath}
\langle \phi,F\mu\rangle_{\Val^\infty(V)}=\langle \phi,\mu\rangle_{\mathcal V^\infty(V)}. 
\end{displaymath}

Fix a Euclidean structure on $V$, which induces canonical identifications $\Dens(V) \cong \mathbb C$ and
$S^*V \cong V \times S^{n-1}$. We may assume by linearity that $\phi$ is $k$-homogeneous. Represent
$\mu=\nu(\omega,\gamma)$ with some compactly supported forms $\omega \in \Omega_c^{n-1}(S^*V)$, $\gamma
\in \Omega_c^n(V)$. 

We may write
\begin{displaymath}
F\mu=\nu\left(\int_V x^* \omega d\vol(x),\int_V x^*\gamma d\vol(x)\right).
\end{displaymath}

If $k<n$, then $\phi=\nu(\beta,0)$ for some form $\beta\in\Omega^{n-1}(S^*V)^{tr}$
by the irreducibility theorem \cite{alesker_mcullenconj01}. By the product formula from \cite{alesker_bernig} we have  
\begin{equation}
\langle \phi,\mu\rangle_{\mathcal V^{\infty}(V)}=\int_{S^*V}\omega\wedge s^{*}D\beta+\int_{V}\gamma\wedge\pi_{*}\beta\label{eq:manifolds_pairing}
\end{equation}
and 
\begin{equation}
\langle \phi, F\mu\rangle_{\Val^{\infty}(V)}=\pi_{*}\left(\int_{V}x^{*}\omega d\vol(x)\wedge s^{*}D\beta\right)
+\int_{S^{n-1}}\beta\cdot\int_V x^{*}\gamma d\vol(x). \label{eq:valuation_pairing}
\end{equation}

The second summand in \eqref{eq:manifolds_pairing} is
\begin{displaymath}
\int_{V}\gamma\wedge\pi_{*}\beta=\int_{S^{n-1}} \beta \cdot \int_V \gamma
\end{displaymath}
which coincides with the second summand of \eqref{eq:valuation_pairing}. 

Denoting $\psi:=s^{*}D\beta \in\Omega^{n}(S^*V)^{tr}$ and $\tau:=\omega \wedge \psi \in \Omega_{c}^{2n-1}(S^*V)$, it
remains to verify that 
\begin{displaymath}
\pi_{*}\left(\int_{V}x^{*}\omega d\vol(x) \wedge\psi\right)=\int_{S^*V}\omega\wedge\psi,
\end{displaymath}
which is equivalent to 
\begin{displaymath}
\pi_{*}\left(\int_{V} x^{*} \tau d\vol(x)\right)=\left(\int_{S^*V} \tau \right)\vol.
\end{displaymath}

Write $\tau=f(y,\theta)d\vol(y)d\theta$ for $y\in V$, $\theta\in S^{n-1}$ and $d\theta$ the volume form on
$S^{n-1}$. Then
\begin{align*}
\pi_{*}\left(\int_{V}x^{*}\tau
d\vol(x)\right)(y) & =\left(\int_{S^{n-1}} \left(\int_{V}f(y+x,\theta)d\vol(x)\right) d\theta \right)\vol(y)\\
& =\left(\int_{S^{n-1}}\left(\int_{V}f(x,\theta)d\vol(x)\right) d\theta \right)\vol(y)\\
& = \left(\int_{S^*V} \tau \right)\vol(y),
\end{align*}
as required.

Now assume $k=n$, so $\phi=\nu(0,\lambda \vol)$ is a Lebesgue measure on $V$. Then 
\begin{displaymath}
\langle \phi, F\mu\rangle_{\Val^{\infty}(V)}=\lambda \pi_*\left(\int_V x^*\omega d\vol(x)\right)
\end{displaymath}
and 
\begin{displaymath}
\langle \phi,\mu\rangle_{\mathcal V^{\infty}(V)}=\lambda \int_{V} \pi_{*}\omega \vol.
\end{displaymath}
Since the right hand sides coincide, the commutativity of the diagram follows.


We proceed to show surjectivity of $F$. Let $\phi$ be a smooth translation invariant valuation on $V$. We fix translation invariant differential
forms $\omega \in \Omega^{n-1}(S^*V), \gamma \in \Omega^n(V)$ with $\phi=\nu(\omega,\gamma)$. 

Let $\vol$ be a density on $V$ and let $\beta$ be a compactly supported smooth function
such that $\int_V \beta(x) d \vol(x)=1$. 

The valuation $\mu:=\nu(\pi^*(\beta) \wedge \omega,\beta \gamma)$ is smooth, compactly
supported and
satisfies $F(\mu)=\phi \otimes  \vol^*$. This shows that $F$ is onto. Thus $F$ induces an isomorphism  
\begin{displaymath}
 \tilde F: \mathcal{V}_c^\infty(V)/\ker F \stackrel{\cong}{\longrightarrow} \Val^\infty(V) \otimes \Dens(V^*).
\end{displaymath}
The transpose of $\tilde F$ is an isomorphism 
\begin{displaymath}
 \tilde F^*: \Val^{-\infty}(V) \to (\ker F)^\perp.
\end{displaymath}
The proof will be finished once we can show $(\ker F)^\perp = \mathcal{V}^{-\infty}(V)^{tr}$.

If $\mu \in \mathcal{V}_c^\infty(V)$, then $(t_v)_*\mu-\mu \in \ker F$ for every $v \in V$, where $t_v$
is the translation by $v$. It follows that $(\ker F)^\perp \subset \mathcal{V}^{-\infty}(V)^{tr}$.

Let $\phi \in \mathcal{V}^{-\infty}(V)^{tr}$. Fix a compactly supported approximate identity $f_\epsilon$ in $\GL(n)$
and
set $\phi_\epsilon:=\phi * f_\epsilon$. Then $T(\phi_\epsilon)=T(\phi) * f_\epsilon$ and
$C(\phi_\epsilon)=C(\phi)*f_\epsilon$ are smooth currents. By \cite[Lemma 8.1]{alesker_bernig}, $\phi_\epsilon \in
\Val^\infty(V)$ and $\phi_\epsilon \to \phi$. This
shows that $\Val^\infty(V) \cong \mathcal{V}^\infty(V)^{tr}$ is dense in $\mathcal{V}^{-\infty}(V)^{tr}$. 

Clearly $\Image\left(F^*:\Val^\infty(V) \to \mathcal{V}^{-\infty}(V)^{tr}\right) \subset (\ker F)^\perp$. Since the
image
is dense in $\mathcal{V}^{-\infty}(V)^{tr}$, it follows that $\mathcal{V}^{-\infty}(V)^{tr} \subset (\ker F)^\perp$.
\endproof  

\begin{Definition}
For $\phi \in \Val^{-\infty}(V)$ we set $\WF(\phi):=\WF(T(\phi))
\subset S^*(S^*V)$. Given $\Gamma \subset S^*(S^*V)$ a closed set, we define 
\begin{displaymath}
 \Val^{-\infty}_\Gamma(V):=\left\{\phi \in \Val^{-\infty}(V): \WF(\phi) \subset \Gamma\right\}.
\end{displaymath}
\end{Definition}

\begin{Lemma} \label{lemma_density_smooth}
 The subspace $\Val^\infty(V) \subset \Val^{-\infty}_\Gamma(V)$ is dense.  
\end{Lemma}

\proof
 This follows by \cite[Lemma 8.2]{alesker_bernig}, noting that in the proof of that lemma, a translation
invariant generalized valuation is approximated by translation invariant smooth valuations. 
\endproof

\section{Embedding McMullen's polytope algebra}
\label{sec_embedding_polytopes}

Let $V$ denote an $n$-dimensional real vector space, and $\Pi(V)$
the McMullen polytope algebra on $V$. It is defined as the abelian group generated by
polytopes,
with the relations of the inclusion-exclusion principle and translation invariance. The product is defined on
generators by $[P] \cdot [Q]:=[P+Q]$. It is almost a graded algebra over $\mathbb{R}$. We refer to
\cite{mcmullen_polytope_algebra} for a detailed study of its properties. 

For $\lambda \in \R$, the dilatation $\Delta(\lambda)$ is defined by $\Delta(\lambda)[P]=[\lambda P]$. The {\it $k$-th
weight space} is defined by  
\begin{displaymath}
 \Xi_k:=\{x \in \Pi: \Delta(\lambda)x=\lambda^k x \text{ for some rational } \lambda>0, \lambda \neq 1\}.
\end{displaymath}

McMullen has shown (\cite{mcmullen_polytope_algebra}, Lemma 20) that 
\begin{displaymath}
 \Pi(V)=\bigoplus_{k=0}^n \Xi_k.
\end{displaymath}

The aim of this section is to prove Theorem \ref{mainthm_injection_mcmullen}, namely that $\Pi(V)$ embeds into
$\Val^{-\infty}(V) \otimes \Dens(V^*)$.

Recall that the map $M:\Pi(V) \to \Val^{-\infty}(V) \otimes \Dens(V^*) \cong  \Val^{\infty}(V)^*$ is
defined by 
\begin{displaymath}
\langle M([P]),\phi\rangle=\phi(P),  \quad  P \in \mathcal{P}(V).
\end{displaymath}

We will denote by $M_k:\Pi(V)\to \Val_k(V) \otimes \Dens(V^*)$ the $k$-homogeneous component of
the image of $M$, and similarly $[\phi]_{k}$ is the $k$-homogeneous component
of a valuation $\phi$.

\begin{Lemma}
\begin{displaymath}
M([P])=\vol(\bullet-P)\otimes \vol^*.
\end{displaymath}
In particular $\Image(M) \subset \Val(V)\otimes \Dens(V^*)$. 
\end{Lemma}

\proof 
We claim that $\langle \PD(\psi),\vol(\bullet - K) \otimes \vol^* \rangle=\psi(K)$ for $\psi
\in \Val^\infty(V)$ and $K \in \mathcal{K}(V)$, where $\PD:\Val^\infty(V) \to \Val^{-\infty}(V)$ denotes the natural embedding given by the Poincar\'e duality of Alesker \cite{alesker04_product}. 

Indeed, let first $K$ be smooth with positive curvature. Then, by (\cite{bernig_aig10}, (15)), 
\begin{displaymath}
 \psi \cdot \vol(\bullet-K)= \int_V \psi((y+K) \cap \bullet) d\vol(y) 
\end{displaymath}
and hence 
\begin{displaymath}
 \langle \PD(\psi), \vol(\bullet-K) \otimes \vol^*\rangle=\psi(K).
\end{displaymath}
By approximation, this holds for non-smooth $K$ as well, showing the claim. Since by definition
$\langle M[P],\psi\rangle=\psi(P)$,
the
statement follows.
\endproof

\begin{Lemma} \label{lemma_m_as_average}
 Let $P$ be a polytope and $\Gamma(P) \in \mathcal{V}^{-\infty}(V)$ the associated generalized valuation, i.e.
$\langle \Gamma(P), \phi\rangle=\phi(P)$ for $\phi \in \mathcal{V}^\infty_c(V)$. Then 
\begin{displaymath}
 F^*M([P])=\int_V \Gamma(P+x) d \vol(x) \otimes \vol^*,
\end{displaymath}
with $F$ as in Proposition \ref{prop_identification_vals}.
\end{Lemma}

\proof
Let $\phi \in \mathcal{V}_c^\infty(V)$. Then 
\begin{align*}
 \langle F^* M([P]),\phi\rangle & =\langle M([P]),F\phi\rangle\\
& = F(\phi)(P)\\
& = \int_V \phi(P+x)d\vol(x) \otimes  \vol^*\\
& = \int_V \langle \Gamma(P+x),\phi\rangle d\vol(x) \otimes \vol^*\\
& = \left\langle \int_V \Gamma(P+x) d\vol(x) \otimes \vol^*,\phi \right\rangle.
\end{align*}
\endproof

For the following, we fix a Euclidean structure and an orientation on $V$ and denote by $\vol$ the corresponding Lebesgue measure on
$V$.

Denote by $\mathcal{C}_{j}$ the collection of $j$-dimensional oriented submanifolds $N \subset S^{n-1}$,
obtained by intersecting $S^{n-1}$ with a $(j+1)$-dimensional polytopal cone $\hat{N}$ in $V$.

 Given $v\in  \Lambda^k V$, define the current $[v] \in \mathcal{D}_{k}(V)$ by 
\begin{displaymath}
 \langle [v],\omega\rangle:=\int_V \omega|_x(v) d\vol(x), \omega \in \Omega_c^k(V).
\end{displaymath}

Let $\Lambda_s^kV$ denote the cone of simple $k$-vectors in $V$. Given a pair $(v,N)
\in \Lambda_s^k V \times \mathcal{C}_{n-k-1}$, we define the current $A_{v,N} = [v]
 \times [[N]]\in
\mathcal{D}_{n-1}(V \times S^{n-1})$, where $[[N]]$ is the current of integration over $N$.

Observe that changing the
sign of $v$ and the orientation of $N$ simultaneously leaves the current $A_{v,N}$ invariant.

Let $Y_k \subset \mathcal{D}_{n-1}(V \times S^{n-1})$ be the $\mathbb{C}$-span of currents of the form
$A_{v,N}, v \in \Lambda^k_sV, N \in \mathcal{C}_{n-k-1}$ such that $\Span(v)\oplus
\Span(\hat{N})=V$ as oriented spaces. 

Given a polytope $P \subset V$, we let $\mathcal{F}_k$ be the set of $k$-faces of $P$. Each face is assumed to possess some fixed orientation. If $F$ is a face of $P$ of
dimension strictly less than $n$, we let $n(F,P)$ be the normal cone of $F$, and $\check n(F,P):=n(F,P) \cap S^{n-1}$ is oriented so that the linear space parallel to $F$, followed by $\Span(\check n(F,P))$, is positively oriented.
Moreover, let $v_F$ be the unique $k$-vector in the linear space parallel to $F$ such that $|v_F|=\vol_k
F$, the sign determined by the orientation of $F$.

\begin{Lemma} \label{lemma_image_of_polytope}
Let $P$ be a polytope. Then  
\begin{align*}
E(M_0([P])) & =\left(0,\vol(P)  [[V]]\right),\\
E(M_{n-k}([P])) & =\left(\sum_{F \in \mathcal{F}_k(P)} A_{v_{F},\check{n}(F,P)},0\right), \quad 0\leq k\leq
n-1. 
\end{align*}
\end{Lemma}

\proof

It follows from \eqref{eq_currents_for_polytope} and Lemma \ref{lemma_m_as_average} that 
\begin{align*}
T(M([P])) & = \int_{V} T(\Gamma(P+x))d\vol(x)=\sum_{F \in \mathcal{F}_{\leq n-1}(P)} A_{v_F,\check{n}(F,P)}\\
C(M([P])) & = \int_{V} C(\Gamma(P+x))d\vol(x) = \vol(P) [[V]].
\end{align*}
From this the statement follows. 
\endproof 

Let 
\begin{displaymath}
T_k:\Image(M_{n-k}) \to Y_k, \quad 0\leq k\leq n-1
\end{displaymath}
be the restriction of $T$ to $\Image(M_{n-k})$. For a linear subspace
$L\subset V$, we will write $S(L)=S^{n-1} \cap L$.

 Let $\mathcal{F}(V)$ denote the space  of $\mathbb{Z}$-valued constructible functions on $V$, i.e. functions of
the form
$\sum_{i=1}^N n_i 1_{P_i}$ with $n_i \in \mathbb{Z}, P_i \in \mathcal{P}(V)$  (compare \cite{alesker_val_man4}).
Let $\mathcal{F}_\text{a.e.}(V)$ be the
set of congruence classes of constructible functions where $f \sim g$ if $f-g=0$ almost everywhere.  

\begin{Lemma} \label{lemma_inclusion_into_constr_functions}
Denote by $Z$ the abelian group generated by all formal integral combination of compact convex polytopes in
$V$. Let $W \subset Z$ denote the subgroup generated by lower-dimensional polytopes and elements of the form 
$[P \cup
Q]+[P \cap Q]-[P]-[Q]$ where
$P \cup Q$ is convex. Then the map 
\begin{align*}
Z/W & \to \mathcal{F}_\text{a.e.}(V)\\
\sum_i n_i [P_i] & \mapsto \sum_i n_i 1_{P_i}, \quad n_i \in \mathbb{Z} 
\end{align*}
 is injective.
\end{Lemma}

\proof

It is easily checked that the map is well-defined. To prove injectivity, it is enough to prove that $f:=\sum_i 
1_{P_i} \sim \sum_j 1_{Q_j}$ implies $\sum_i [P_i] \equiv \sum_j [Q_j]$. 
Decompose the connected components of 
\begin{displaymath}
\left(\cup_i P_i \cup \cup_j Q_j\right) \setminus \left(\cup_i \partial P_i \cup \cup_j \partial Q_j\right) 
\end{displaymath}
into simplices $\{\Delta\}$, disjoint except at their boundary. Then, by the inclusion-exclusion principle, 
\begin{displaymath}
\sum_i [P_i] \equiv \sum_i \sum_{\Delta \subset P_i} [\Delta] \equiv \sum_{k \geq 1} (-1)^{k+1} \sum_{\Delta \subset \bigcup_{i_1<\ldots<i_k} \cap_{j=1}^k P_{i_j}}[\Delta]
\end{displaymath}

By examining the superlevel set $\{f \geq k\}$ we see that 
\begin{displaymath}
\bigcup_{i_1<\ldots<i_k}(P_{i_1} \cap \ldots \cap P_{i_k})  \sim \bigcup_{j_1<\ldots<j_k}(Q_{j_1} \cap \ldots \cap
Q_{j_k}),
\end{displaymath}
where $A \sim B$ means the sets $A,B$ coincide up to a set of measure zero.

Therefore, 
\begin{displaymath}
 \sum_i [P_i] \equiv \sum_j [Q_j] \mod W,
\end{displaymath}
as claimed.
\endproof

\begin{Remark}
The same claim and proof apply if we replace $Z$ with the free abelian
group of polytopal cones with vertex in the origin.
\end{Remark}

Let us recall some notions from \cite{mcmullen_polytope_algebra}. Let $L$ be a subspace of $V$. The {\it cone group}
$\hat \Sigma(L)$ is the
abelian group with generators $[C]$, where $C$ ranges over all convex polyhedral cones in $L$, and with the relations
\begin{enumerate}
 \item $[C_1 \cup C_2]+[C_1 \cap C_2]=[C_1]+[C_2]$ whenever $C_1,C_2,C_1 \cup C_2$ are convex polyhedral cones;
\item $[C]=0$ if
$\dim C < \dim L$. 
\end{enumerate}

The {\it full cone group} is given by 
\begin{displaymath}
 \hat \Sigma:=\bigoplus_{L \subset V} \hat \Sigma(L),
\end{displaymath}
where the sum extends over all linear subspaces of $L$. 

\begin{Lemma}[Lemma 39 from \cite{mcmullen_polytope_algebra}] \label{lemma_lemma39}
 The map 
\begin{align*}
 \sigma_k : \Pi(V) & \to \mathbb{C} \otimes_ \mathbb{Z} \hat \Sigma \\
[P] & \mapsto \sum_{F \in \mathcal{F}_k(P)} \vol(F) \otimes n(F,P)
\end{align*}
restricts to an injection on $\Xi_k$. 
\end{Lemma}

\begin{Proposition}
For all $0\leq k\leq n-1$, there is a linear map $\Phi_k:Y_k \to\mathbb{C} \otimes_{\mathbb{Z}} \hat \Sigma$
such that $\Phi_k(A_{v,N})=|v| \otimes [\hat{N}]$.
\end{Proposition}

\proof
The first thing to note is that if $A_{v_1,N_1}=A_{v_2,N_2}$, then either $v_1=v_2$,
$N_1=N_2$ or $v_1=-v_2$ and $N_1=\overline{N_2}$, where
$\overline{N_j}$ is $N_j$ with reversed orientation. Thus on the generators
of $Y_k$, $\Phi_k(A_{v,N}):=|v| \otimes [\hat{N}]$ is well-defined. 

Now assume that $\sum_j c_j A_{v_j,N_j}=0$. We shall show that $\sum_j c_j |v_j| \otimes [\hat{N}_j]=0$.

Note that for all $\rho \in \Omega_c^k(V)$ and $\omega \in \Omega^{n-k-1}(S^{n-1})$, one has 
\begin{displaymath}
\sum_j c_j \int_V \rho|_x(v_j)d\vol(x) \int_{N_j} \omega=0.
\end{displaymath}

More generally, suppose that we have 
\begin{displaymath}
\sum_j \lambda_j \int_{N_j} \omega=0
\end{displaymath}
for some coefficients $\lambda_j \in \mathbb{C}$ and for all $\omega\in\Omega^{n-k-1}(S^{n-1})$. 

Let $L \subset V$ be a linear subspace of dimension $n-k$. Let $U_\epsilon \subset S^{n-1}$ denote the
$\epsilon$-neighborhood of $L \cap S^{n-1}$, where $\epsilon$ is sufficiently small. Let $p_\epsilon:U_\epsilon \to
L\cap S^{n-1}$ denote the nearest-point projection. 

Fix some $\beta \in C^{\infty}[0,1]$ such that $\beta(x)=1$ for $0 \leq x \leq \frac13$ and $\beta(x)=0$
for $\frac23 \leq x \leq 1$. Let $ \beta_\epsilon \in C^{\infty}(U_\epsilon)$ be
given by $\beta_\epsilon(x)=\beta(\dist(x,p(x))/\epsilon)$. 

Given a form $\sigma \in \Omega^{n-k-1}(L\cap S^{n-1})$, we set
$\omega:=\beta_\epsilon p_\epsilon^* \sigma \in \Omega^{n-k-1}(S^{n-1})$. 

If $N_j \subset L \cap S^{n-1}$, then   
\begin{displaymath}
\int_{N_j} \omega=\int_{N_j} \sigma,
\end{displaymath}
while if $N_j$ does not lie in $L \cap S^{n-1}$ then $\int_{N_j} \omega \to 0$
as $\epsilon \to 0$. 

Thus, letting $\epsilon \to 0$, we obtain
\begin{displaymath}
\sum_{N_j \subset L} \lambda_j \int_{N_j} \sigma=0,
\end{displaymath}

where the sum is over all $N_j$ contained in $L$.

Now fix an orientation on
$L$ and assume without loss of generality that all $N_j \subset L$ have the induced orientation. Going back to the
original equation we may write 
\begin{displaymath}
\sum_{N_j \subset L} c_j \int_V \rho|_x(v_j)d\vol(x) \int_{N_j}\sigma=0.
\end{displaymath}

Let $v_0 \in \Lambda^k L^\perp \subset \Lambda^k V$ be the unique vector with $|v_0|=1$ and $v_0^\perp=L$ (the last
equality understood with orientation). Then $v_j=|v_j|v_0$ for all $j$  with $N_j \subset L$. Therefore, the
above equation implies that
\begin{displaymath}
\sum_{N_j \subset L} c_j |v_j| \int_{N_j}\sigma=0
\end{displaymath}
for all $\sigma \in \Omega^{n-k-1}(L \cap S^{n-1})$ and this implies that 
\begin{displaymath}
\sum_{N_j \subset L} c_j |v_j| 1_{N_j}=0
\end{displaymath}
almost everywhere. 

By Lemma \ref{lemma_inclusion_into_constr_functions} we have in $\mathbb{C} \otimes \hat{\Sigma}(L)$ 
\begin{displaymath}
\sum_{N_j \subset L} c_j |v_j| \otimes [\hat{N}_j]=0.
\end{displaymath}
Since $L$ was arbitrary, we deduce that in $\mathbb{C} \otimes \hat{\Sigma}$ we have 
\begin{displaymath}
\sum_j c_j |v_j| \otimes [\hat N_j]=0,
\end{displaymath}
as required.
\endproof

\proof[Proof of Theorem \ref{mainthm_injection_mcmullen}]
By Lemma \ref{lemma_lemma39}, the map 
\begin{displaymath}
\sigma_k:\Xi_k \to \mathbb{C} \otimes \hat{\Sigma}
\end{displaymath}
is injective. For $0 \leq k \leq n-1$ we have   
\begin{displaymath}
\sigma_k=\Phi_k \circ T_k \circ M_{n-k},
\end{displaymath}
while $\sigma_n$ is the volume functional on $\Pi(V)$ restricted to
$\Xi_n$, so identifying the space of Lebesgue measures on $V$ with $\mathbb{C}$ allows us to write $\sigma_n=C \circ
M_0$. We conclude that $(M_{n-k})|_{\Xi_k}$ is injective for each $k$, and hence $M$ is injective.
\endproof

\section{Partial convolution product}
\label{sec_partial_conv}

Let us first describe the convolution in the smooth case, see \cite{bernig_fu06}, but using a more intrinsic approach. 

Let $V$ be an $n$-dimensional vector space. Let $\Dens(V)$ denote the $1$-dimensional $\GL(V)$-module
 of densities on $V$. The orientation bundle $\ori(V)$ is the real $1$-dimensional linear space consisting of all functions $\rho:\Lambda^{top}V \to \mathbb R$ such that $\rho (\lambda \omega)=\mathrm{sign}(\lambda) \rho(\omega)$ for all $\lambda \in\mathbb R$ and $\omega \in \Lambda^{top}V$. Note that there is a canonical isomorphism $\ori(V) \cong \ori(V^*)$.

Let
$\mathbb P_+(V)$ be the space of oriented $1$-dimensional subspaces of $V$. Then $ S^*V:=V \times \mathbb P_+(V^*)$
has a natural contact
structure. 

We have a natural
non-degenerate pairing  
\begin{displaymath}
\Lambda^k V^* \otimes \Lambda^{n-k}V^* \stackrel{\wedge}{\longrightarrow} \Lambda^n V^* \cong
\Dens(V) \otimes \ori(V),
\end{displaymath}
which induces an isomorphism
\begin{displaymath}
 *:\Lambda^k V^* \otimes \ori(V) \otimes \Dens(V^*) \stackrel{\cong}{\longrightarrow} (\Lambda^{n-k} V^*)^*
\cong \Lambda^{n-k} V. 
\end{displaymath}

Let 
\begin{displaymath}
\ast_1:\Omega (S^*V)^{tr} \otimes \ori(V) \otimes \Dens(V^*) 
\stackrel{\cong}{\longrightarrow} \Omega(V^* \times \mathbb P_+(V^*))^{tr}  
\end{displaymath}
be defined by
\begin{displaymath}
 \ast_1(\pi_1^*\gamma_1 \wedge \pi_2^*\gamma_2):=(-1)^{\binom{n-\deg \gamma_1}{2}} \pi_1^*(*\gamma_1) \wedge
\pi_2^*\gamma_2
\end{displaymath} 
for $\gamma_1\in \Omega(V)^{tr} \otimes \ori(V) \otimes \Dens(V^*), \gamma_2 \in \Omega(\mathbb
P_+(V^*))$.

Let $\phi_j=\nu(\omega_j,\gamma_j) \in
\Val^\infty(V) \otimes \Dens(V^*),j=1,2$, where $\omega_j \in \Omega^{n-1}(S^*V)^{tr} \otimes
\ori(V) \otimes \Dens(V^*), \gamma_j \in \Omega^n(V)^{tr} \otimes \ori(V) \otimes \Dens(V^*)$. Then the convolution
product
$\phi_1 * \phi_2$ is defined as $\nu(\omega,\gamma) \in \Val^\infty(V) \otimes \Dens(V^*)$, where $\omega$ and $\gamma$
satisfy 
\begin{align*}
 D\omega +\pi^*\gamma & = \ast_1^{-1}\left(\ast_1 (D\omega_1+\pi^*\gamma_1) \wedge
\ast_1(D\omega_2+\pi^*\gamma_2)\right)\\
 \pi_*\omega & = \pi_* \circ \ast_1^{-1} \left( \ast_1 \kappa_1 \wedge \ast_1 (D\omega_2)\right)+(*\gamma_1)\pi_*\omega_2+(*\gamma_2)\pi_*\omega_1.
\end{align*}
Here $\kappa_1 \in \Omega^{n-1}(S^*V)^{tr}$ is any form such that $d\kappa_1=D\omega_1$. 

The convolution extends to a partially defined convolution on the space $\Val^{-\infty}(V) \otimes \Dens(V^*)$ as
follows. 

The space $\mathcal{D}_*(S^*V)$ admits a bigrading
\begin{displaymath}
\mathcal{D}_*(S^*V)=\bigoplus_{k=0}^n \bigoplus_{l=0}^{n-1} \mathcal{D}_{k,l}(S^*V),
\end{displaymath}
and for $T \in\mathcal{D}_*(S^*V)$, we denote by $[T]_{k,l}$
the component of bidegree $(k,l)$.

We consider now the $\GL(V)$-module of translation-invariant currents $\mathcal D(S^*V)^{tr}=\mathcal
D(V)^{tr} \otimes \mathcal D(\mathbb P_+(V^*))$. One has the natural identification $\mathcal
D_k(V)^{tr}=\Omega^{n-k}(V)^{tr} \otimes \ori(V)$, and a non-degenerate pairing $\mathcal D_k(V)^{tr} \otimes
\Omega^k(V)^{tr}\to \Dens(V)$.

We define for $T \in \mathcal{D}_{k,l}(S^*V)^{tr} \otimes \ori(V) \otimes \Dens(V^*)$ the element
$\ast_1 T \in \mathcal{D}_{n-k,l}(V^*\times \mathbb P_+(V^*))^{tr}$ by
\begin{displaymath}
\langle \ast_1 T,\delta\rangle:=(-1)^{nk+nl+k+{\binom{n}{2}}} \langle T,\ast_{1}^{-1}\delta\rangle \in \Dens(V^*)
\end{displaymath}
for all $\delta \in \Omega_c^{n-k,l}(V^*\times \mathbb P_+(V^*))^{tr}$. With this definition of $\ast_1$ the
diagram 
\begin{displaymath}
\xymatrixcolsep{2pc}
 \xymatrix{\Omega^{k,l}(S^*V)^{tr} \otimes \ori(V) \otimes \Dens(V^*)  \ar[r]^-{*_1}
\ar@{_{(}->}[d] & 
\Omega^{n-k,l}(V^*\times \mathbb P_+(V^*))^{tr} \ar@{_{(}->}[d] \\
\mathcal{D}_{n-k,n-l-1}(S^*V)^{tr}  \otimes \ori(V)\otimes \Dens(V^*)\ar[r]^-{*_1}&  
\mathcal{D}_{k,n-l-1}(V^*\times \mathbb P_+(V^*))^{tr}}
\end{displaymath}
commutes. Equivalently,
\begin{displaymath}
 \ast_1 \left([[S^*V]] \llcorner \gamma\right)=[[S^*V]] \llcorner \ast_1
\gamma
\end{displaymath}
for all $\gamma \in \Omega(S^*V)^{tr} \otimes \ori(V) \otimes \Dens(V^*)$.
Clearly we have $\WF(\ast_{1} T)=\WF(T)$ for $T\in \mathcal{D}(S^*V)^{tr}$

\begin{Remark}
Without the choice of an orientation and a Euclidean scalar product, we may intrinsically define $A_{v,N} \in \mathcal
D_{k,n-1-k}(S^*V)^{tr} \otimes \ori(V) \otimes \Dens(V^*)$ as follows. 

Let $v \in \Lambda^k_s V$ and let $N \subset \mathbb P_+(V^*)$ be an $(n-1-k)$-dimensional, geodesically
convex polytope contained in $v^\perp \cap \mathbb P_+(V^*)$. 
Then define  
\begin{displaymath}   
 \langle A_{v,N}, \pi_1^*\gamma \wedge \pi_2^* \delta  \otimes \sigma \otimes \epsilon \rangle=
\int_V \langle x^*\omega, v\rangle d\sigma(x)\int_{N_\epsilon} \delta
\end{displaymath}
for $\gamma \in \Omega_c^k(V), \delta \in \Omega^{n-k-1}(\mathbb P_+(V^*)), \epsilon\in\ori(V)$ and $\sigma \in
\Dens(V)$. Here $N_\epsilon$ equals $N$ with a choice of orientation such that the orientation of the pair 
$(v,N_\epsilon)$ is induced by $\epsilon$.

Given a Euclidean trivialization, this reduces to the previous definition of $A_{v,N}$.
\end{Remark}

\begin{Lemma} \label{lemma_filling_current}
Given a Legendrian cycle $T\in\mathcal{D}_{n-k,k-1}(S^*V)^{tr}$
with $1\leq k\leq n-1$, there exists $\tilde{T} \in \mathcal{D}_{n-k,k}(S^*V)^{tr}$ with  $T=\partial
\tilde{T}$ and
$\WF(\tilde{T})=\WF(T)$.
\end{Lemma}

\proof
Let us use a Euclidean scalar product and an orientation on $V$. Then we may identify $\Dens(V) \cong \mathbb{C},
\ori(V) \cong \mathbb{C}, S^*V \cong SV=V \times S^{n-1}$.

Let $\phi \in \Val_{k}^{-\infty}(V)$ be the valuation represented by $(T,0)$.
Choose a sequence $\phi_j \in \Val_k^\infty(V)$ such that $\phi_j \to \phi$. Let $\omega_j \in
\Omega^{k,n-k-1}(SV)^{tr}=\Omega^{n-k-1}(S^{n-1}) \otimes \Lambda^kV^*$ be a form representing $\phi_j$ and such that
$D\omega_j=d\omega_j$. 

Note that $\mathcal{D}_{n-k,k-1}(SV)^{tr}=\mathcal{D}_{k-1}(S^{n-1}) \otimes \Lambda^{n-k} V$. Then $d\omega_j
\in \Omega^{n-k}(S^{n-1}) \otimes \Lambda^kV^* \subset \mathcal{D}_{k-1}(S^{n-1}) \otimes \Lambda^{n-k}V$
converges weakly to $T$. 

Let  $G:\Omega^*(S^{n-1}) \to \Omega^*(S^{n-1})$ denote the Green operator on $S^{n-1}$ and
$\delta:\Omega^*(S^{n-1}) \to \Omega^{*-1}(S^{n-1})$ the codifferential. We define $\beta_j:=G(d\omega_j) \in
\Omega^{n-k}(S^{n-1}) \otimes
\Lambda^{n-k}V$, i.e. $\Delta \beta_j=d\omega_j$. Then $\Delta d\beta_j=d\Delta \beta_j=0$, hence $d\beta_j$ is
harmonic, which implies that $\delta d\beta_j=0$. 

We define 
\begin{displaymath}
 \tilde T:=(-1)^{n+k}\lim_{j \to \infty}  [[S^{n-1}]] \llcorner \delta\beta_j \in
\mathcal{D}_k(S^{n-1}) \otimes \Lambda^{n-k}V \subset
\mathcal{D}_{n-k,k}(SV)^{tr}.
\end{displaymath}

Then
\begin{displaymath}
\partial \tilde T =  \lim_j  [[S^{n-1}]] \llcorner d \delta\beta_j = \lim_j  [[S^{n-1}]] \llcorner \Delta
\beta_j=\lim_j [[S^{n-1}]] \llcorner d\omega_j =T. 
\end{displaymath}

From $\tilde T=\delta \circ G(T)$ and \eqref{eq_wavefront_laplacian}, we infer that $\WF(\tilde T) \subset
\WF(T)$.
Conversely, from \eqref{eq_wavefront_boundary} we deduce that $\WF(T)=\WF(\partial \tilde{T}) \subset
\WF(\tilde T)$.
\endproof
 
Recall that $ s:S^*V \to S^*V$ is the antipodal map. 

\begin{Definition}
Let $\phi_j \in \Val^{-\infty}(V)  \otimes \Dens(V^*), E(\phi_j)=:(T_j,C_j)$, $j=1,2$. We call $\phi_1,\phi_2$
transversal if 
\begin{displaymath}
 \WF(T_1) \cap  s(\WF(T_2)) = \emptyset. 
\end{displaymath}
\end{Definition}

\begin{Proposition} \label{prop_def_convolution}
Let $\phi_j \in \Val^{-\infty}(V) \otimes \Dens(V^*)$, $j=1,2$ be transversal and $(T_j,C_j):=E(\phi_j), j=1,2$.
Decompose
$T_j=t_j+T'_j$, where $t_j=\alpha_j  \pi^*([[V]] \llcorner \vol_n) \otimes \vol_n^* \in \mathcal{D}_{0,n-1}(S^*V)^{tr} \otimes \Lambda^n V^*$ is the
corresponding $(0,n-1)$-component, and let $\tilde T_1 \in \mathcal{D}_n(S^*V)^{tr} \otimes \Lambda^n V^*$ be a current such
that $\partial \tilde T_1=T'_1$ and $\WF(\tilde T_1)=\WF(T_1)$, guaranteed to exist by Lemma
\ref{lemma_filling_current}. Then the currents 
\begin{align*}
T & :=\ast_1^{-1}\left(\ast_1T_1 \cap \ast_1T_2\right)\\
C & :=\pi_*\left(\ast_1^{-1}\left(\ast_1\tilde{T}_1 \cap \ast_1
T'_2\right)\right)+\alpha_1C_2+\alpha_2C_1
\end{align*}
are independent of the choice of $\tilde T_1$ and satisfy the conditions \eqref{eq_conditions_t_c}. 
\end{Proposition}

\proof
Note first that $\partial$ commutes (up to sign) with $*_1$ on translation invariant currents. By \eqref{eq_boundary_intersection}, $T$ equals (up to
sign) the boundary of the $n$-current $S:=\ast_1^{-1}\left(\ast_1 \tilde T_1 \cap \ast_1T_2\right)$. In particular, $T$
is a cycle. Moreover, $\pi_*T= \pm \partial \pi_* S=0$, since $\pi_*S$ is a translation invariant $n$-current,
hence a
multiple of the integration current on $V$  which has no boundary.  For the same reason, $\partial
C=0$, hence the
condition $\pi_*T=\partial C$ is trivially satisfied. 

Note that whenever $Q \in \mathcal{D}_n(S^*V)$ is a boundary, then $\pi_*Q=0$. Indeed, let $Q=\partial R$ and $\rho \in
\Omega^n_c(V)$. Then 
\begin{equation} \label{eq_push_forward_boundary}
 \langle \pi_*Q,\rho\rangle=\langle Q,\pi^*\rho\rangle=\langle \partial R,\pi^* \rho\rangle=\langle R,\pi^*
d\rho\rangle=0. 
\end{equation}

By Lemma \ref{lemma_filling_current}, there exists a translation invariant $n$-current $\tilde T_2$ with $\partial
\tilde T_2=T_2'$. Suppose that $\partial \tilde T_1=\partial \hat T_1=T_1'$ for two $n$-currents $\tilde T_1$
and $\hat T_1$. Then the $n$-current $Q:=*_1^{-1} \left(*_1(\tilde T_1-\hat T_1) \cap *_1 T_2'\right)$ is (up to a
sign) the boundary of the $(n+1)$-current $R:=*_1^{-1} \left(*_1(\tilde T_1-\hat T_1) \cap *_1 \tilde T_2\right)$. By
\eqref{eq_push_forward_boundary}, it follows that $\pi_*Q=0$, which shows that $C$ is independent of the choice of
$\tilde T_1$. 

Let us finally show that $T$ is Legendrian. Fix sequences $(\phi_j^i)_i$ of smooth and translation
invariant valuations converging to $\phi_j, j=1,2$. Let $\phi_j^i$ be represented by the forms $(\omega_j^i,
\gamma_j^i)$. Then $E(\phi_j^i)=(T_j^i, C_j^i)$ is given by the formulas
 \eqref{eq_pair_currents_smooth_a}, \eqref{eq_pair_currents_smooth} and hence 
$\phi_1^i * \phi_2^i$ is represented by the current $T^i=[[S^*V]] \llcorner s^*\kappa^i$, with 
\begin{equation} \label{eq_conv_smooth}
s^*\kappa^i:=\ast_1^{-1}(\ast_1 s^*(D\omega_1^i+\pi^*\gamma_1^i) \wedge \ast_1 s^*(D\omega_2^i+\pi^* \gamma_2^i)). 
\end{equation}
It is easily checked that $\kappa^i$ is a horizontal, closed $n$-form (compare also \cite[Eq. (37)]{bernig_fu06}). It
follows that $T^i$ is Legendrian. 

Note that $ [[S^*V]] \llcorner s^*(D\omega_j^i+\pi^*\gamma_j)$ converges to $T_j$. By the definition of the intersection
current,
$[[S^*V]] \llcorner s^*\kappa^i$ converges to $T$ and hence $T$ is Legendrian.   
\endproof

\begin{Definition} \label{def_convolution}
In the same situation, the convolution product
$\phi_1 * \phi_2 \in
 \Val^{-\infty}(V) \otimes \Dens(V^*)$ is defined as $\phi_1 * \phi_2:=E^{-1}(T,C)$. 
\end{Definition}

\begin{Proposition}
If $\phi_1,\phi_2 \in \Val^\infty(V) \otimes \Dens(V^*) \subset \Val^{-\infty}(V) \otimes \Dens(V^*)$, then the
convolution of Definition \ref{def_convolution} coincides with the convolution from \cite{bernig_fu06}.
\end{Proposition}

\proof

For $\phi_j \in
\Val^\infty(V)$
given by the pairs $(\omega_j, \gamma_j)\in \Omega^{n-1}(S^*V)^{tr} \otimes \Dens(V)$, the corresponding
currents are $E(\phi_{j})=([[S^*V]] \llcorner s^*(D\omega_j+\pi^*\gamma_j), [[V]] \llcorner \pi_*\omega_j)$.

We consider two cases:

If $\omega_1=0$ and $\gamma_1=c\cdot \vol$, then $\phi_1=c \cdot \vol$, and $\phi_1 \ast \phi_2=c\phi_2$
by the original definition of convolution. 

By the new definition, $E(\phi_1)=([[S^*V]] \llcorner \pi^* \gamma_1,0)$, and 
\begin{displaymath}
\ast_1 \pi^* \gamma_1=c \in \Omega^0(V\times S^{n-1}). 
\end{displaymath}

Write 
\begin{displaymath}
\phi_1 \ast \phi_2=E^{-1}( [[S^*V]] \llcorner \pi^*\gamma_1,0) \ast E^{-1}(0,C_2)+E^{-1}( [[S^*V]] \llcorner  \pi^*\gamma_1,0) \ast
E^{-1}(T_2,0). 
\end{displaymath}

If $E^{-1}(0,C_2)=\lambda\chi$, by definition the first summand
equals $c\lambda\chi=cE^{-1}(0,C_2)$. The second summand is $c\cdot E^{-1}(T_2,0)$,
and so $\phi_1 \ast \phi_2=c E^{-1}(T_2,C_2)=c \phi_2$, as
required.

In the remaining case, we may assume $\gamma_1=\gamma_2=0$, so
$E(\phi_{j})=( [[S^*V]] \llcorner  s^{*}D\omega_{j}, [[V]] \llcorner \pi_{*}\omega_{j})$. Moreover, we may assume that
$d\omega_{1},d\omega_{2}$ are vertical.

The original
definition of convolution gives $E(\phi_{1}\ast\phi_{2})=( [[S^*V]] \llcorner s^* D\omega, [[V]] \llcorner \pi_* \omega)$ with 
\begin{align*}
\omega & =\ast_1^{-1}\left(\ast_1 \omega_1 \wedge \ast_1 D\omega_2\right)\\
D\omega & =\ast_1^{-1}\left(\ast_1 D\omega_1 \wedge \ast_1 D\omega_2\right).
\end{align*}

Since $[E^{-1}([[S^*V]] \llcorner s^* D\omega_j,0)]_{n}=[\nu(\omega_{j},0)]_n=0$,
it remains to verify that by the new definition, 
\begin{displaymath}
E^{-1}([[S^*V]] \llcorner  s^* D\omega_1,0) \ast E^{-1}( [[S^*V]] \llcorner  s^* D\omega_2,0)=E^{-1}( [[S^*V]] \llcorner s^* D\omega, [[V]] \llcorner \pi_* \omega).
\end{displaymath}
By homogeneity, we may assume that $\deg \omega_1=(k,n-1-k)$ and
$\deg\omega_{2}=(l,n-1-l)$. 

If $k+l<n$ then by dimensional considerations $\omega=0$. 

If $k+l>n$ then $\pi_* \omega=0$, and the new definition of convolution
gives
\begin{displaymath}
E(\phi_1 \ast \phi_2) =\left( [[S^*V]] \llcorner  \ast_1^{-1}\left((\ast_1s^*D\omega_1) \wedge (\ast_1 s^*
D\omega_2)\right),0\right)=\left( [[S^*V]] \llcorner  s^* D\omega,0\right)
\end{displaymath}
as required. 

Finally, if $k+l=n$, then $1 \leq k \leq n-1$, $D\omega=0$, and by
the new definition $T(\phi_1 \ast \phi_2)=0$. Since $T_1= [[S^*V]] \llcorner  s^* d\omega_1$ and
$\omega_1\in\Omega^{n-1}(S^*V)$, using \eqref{eq_boundary_vs_differential} one can take $\tilde{T}_1=(-1)^n [[S^*V]]
\llcorner  s^* \omega_1$. Using \eqref{eq_push_forward_forms_currents}, the fact that 
the operations $s^* $, $\ast_1$ and $[[S^*V]]\llcorner$  commute, while $\pi_* \circ s^*=(-1)^n
\pi_*$, we obtain 
\begin{align*}
C(\phi_1 \ast \phi_2) & =(-1)^n \pi_* \left( [[S^*V]] \llcorner  \ast_1^{-1}(\ast_1 s^* \omega_1 \wedge \ast_1 s^*
d\omega_2)\right)\\
& =(-1)^n \pi_* ([[S^*V]] \llcorner  s^*\omega)\\
& =(-1)^n [[V]] \llcorner \pi_*  s^*\omega\\
& =[[V]] \llcorner \pi_* \omega,
\end{align*} 
completing the verification.
\endproof

\begin{Proposition} \label{prop_continuity}
Let $\Gamma_1,\Gamma_2 \subset  S^*(S^*V)$ be closed sets with $\Gamma_1 \cap  s\Gamma_2
= \emptyset$ and set 
\begin{displaymath}
 \Gamma:=\Gamma_1 \cup \Gamma_2 \cup \left\{(x,[\xi],[\eta_1+\eta_2]):(x,[\xi]) \in S^*V, (x,[\xi],[\eta_1]) \in
\Gamma_1, (x,[\xi],[\eta_2]) \in \Gamma_2 \right\}.
\end{displaymath}
Then the convolution is a (jointly sequentially) continuous map 
\begin{displaymath}
 \Val^{-\infty}_{\Gamma_1}(V) \otimes \Dens(V^*) \times \Val^{-\infty}_{\Gamma_2}(V) \otimes \Dens(V^*)
\to
\Val^{-\infty}_\Gamma(V)
\otimes \Dens(V^*). 
\end{displaymath}
\end{Proposition}

\proof
In the notations of Proposition \ref{prop_def_convolution}, we have $\WF(*_1 \tilde T_1)=\WF(\tilde T_1) \subset
\WF(T_1) \subset \Gamma_1$ and $\WF(*_1 T_2')=\WF(T_2') \subset \WF(T_2) \subset \Gamma_2$. Since the intersection of
currents is a jointly sequentially continuous map $\mathcal{D}_{*,\Gamma_1}(S^*V) \times \mathcal{D}_{*,\Gamma_2}(S^*V) \to
\mathcal{D}_{*,\Gamma}(S^*V)$, the statement follows.
\endproof

\begin{Corollary}
Whenever it is defined, the convolution is commutative and associative.  
\end{Corollary}

\proof
Let $\phi_j \in \Val^{-\infty}_{\Gamma_i}(V) \otimes \Dens(V^*), j=1,2$. By Lemma \ref{lemma_density_smooth} there
exist sequences $\phi_j^i \in
\Val^\infty(V) \otimes \Dens(V^*), j=1,2$ converging to $\phi_j$ in $\Val^{-\infty}_{\Gamma_j}(V) \otimes \Dens(V^*)$.
By
Proposition \ref{prop_continuity},
$\phi_1^i * \phi_2^i$ converges to $\phi_1 * \phi_2$ in $\Val^{-\infty}_{\Gamma}(V)$, while $\phi_2^i *
\phi_1^i$ converges to $\phi_2 * \phi_1$. Since
the convolution on smooth valuations is commutative, it follows that $\phi_1 * \phi_2=\phi_2 * \phi_1$.

For associativity, let $\Gamma_1,\Gamma_2,\Gamma_3 \subset S^*(S^*V)$ be closed sets such that if $(x,[\xi])
\in S^*V, (x,[\xi],[\eta_i]) \in \Gamma_i, i=1,2,3$, then  
\begin{displaymath}
 \eta_1+\eta_2 \neq 0, \eta_1+\eta_3 \neq 0, \eta_2+\eta_3 \neq 0, \eta_1+\eta_2+\eta_3 \neq 0. 
\end{displaymath}
One easily checks that both maps $\Val^{-\infty}_{\Gamma_1}(V) \times \Val^{-\infty}_{\Gamma_2}(V)  \times
\Val^{-\infty}_{\Gamma_3}(V) \to \Val^{-\infty}(V)$ are well-defined. An approximation argument as above shows
that they agree. 

\endproof

\section{The volume current on the sphere}\label{sec:Appendix}

The aim of the section is the construction of a certain family of currents on the sphere, which can be used to compute the volume of the convex hull of two polytopes on the sphere. This can be viewed as a generalization of the Gauss formula for area in the plane. 
The construction in this section uses tools from geometric measure theory, and is independent of the rest of the paper. It is used in the next section for the proof of the second main theorem.

A geodesically convex polytope on $S^{n-1}$ is the intersection of a proper convex closed polyhedral cone in $V$ with
$S^{n-1}$. We let
$\mathcal{P}(S^{n-1})$ denote the set of oriented geodesically convex polytopes on $S^{n-1}$. For $I \in
\mathcal{P}(S^{n-1})$ of dimension $k$, denote by $E_I=\Span_V I\cap S^{n-1}$ the
$k$-dimensional equator
that it spans. 
For $I,J \subset S^{n-1}$ geodesically convex polytopes, we let $\conv(I,J)
\subset S^{n-1}$ denote the union of all shortest geodesic intervals having an endpoint in $I$ and an endpoint in $J$.
If $\dim
I+\dim J=n-2$, both $I,J$ are oriented and $E_I \cap E_J=\emptyset$, one has a natural orientation on $\conv(I,J)$, by
comparing the orientation of $\Span_V(I) \oplus \Span_V(J)=V$ with the orientation of $V$. The
geodesically convex polytope $-I$ is oriented in such a way that the antipodal map $I \mapsto -I$ is orientation
preserving. Note that $\partial(-I)=-\partial I$ whenever $\dim I>0$, while when $\dim I=0$ we have $\partial
(-I)=\partial I$. Here and in the following, $\partial$ denotes the extended singular boundary operator, i.e.
for a positively oriented point $I$ we have $\partial I=1$. If
$\epsilon=\pm 1$, we write $\conv(I,\epsilon)=I^\epsilon$, where $I^{-1}$ denotes orientation reversal.

We denote $A_{n-1}(I,J):=\vol_{n-1}(\conv(I,J))$ provided that $(-I)\cap J=\emptyset$. Note that whenever the orientation of $\conv(I,J)$ is not well-defined, it is a set of volume zero. Note also that $A_{n-1}$ is a partially-defined bi-valuation in $I,J$, so we may extend $A_{n-1}$ as a partially-defined bilinear functional on chains of polytopes. 

\begin{Lemma} \label{lemma_cancellation}
Let $J \subset S^{n-1}$ be an oriented geodesically convex polytope of dimension $k$. Suppose it does not intersect
$\Span_{\R^n}(I) \cap S^{n-1}$ for all $I \in \partial n_{F}$, for all $F \in \mathcal{F}_k(P)$. Let $\omega \in \Omega_c^k(\mathbb{R}^{n})$. Then 
\begin{equation} \label{eq:conjecture}
\sum_{F \in \mathcal{F}_k(P)} \langle \ [v_F],\omega\rangle A_{n-1}(-\partial
n_F,J)=0. 
\end{equation}
\end{Lemma}

Remark: Note that if $-\partial n_F=\sum_j I_j$ is the decomposition of $\partial n_F$ into geodesically convex polytopes,
then by definition 
\begin{displaymath}
A_{n-1}(-\partial n_F,J)=\sum A_{n-1}(I_{j},J)
\end{displaymath}
is the sum of the oriented volumes. Note also that in formula (\ref{eq:conjecture}) the orientation of $F$ in each summand appears twice, and so the summands
are well-defined.

\proof
For $k=n-1$, this is the well-known statement that 
\begin{displaymath}
\sum_{F\in\mathcal{F}_{n-1}(P)} [v_F]=0,
\end{displaymath}
where the orientation is given by fixing the outer normals to $P$. 

Now for $k<n-1$, fix an arbitrary orientation for all faces of dimension
$k$ and $(k+1)$. For a pair of faces $F \in \mathcal{F}_k(P)$, $G \in \mathcal{F}_{k+1}(P)$, s.t. $F \subset \partial G$, define $\mathrm{sign}(F,G)=\pm 1$
according to the orientation. Then 
\begin{displaymath}
A_{n-1}(-\partial n_F,J)=\sum_{G\supset F} A_{n-1}(-n_G,J)\mathrm{sign}(F,G),
\end{displaymath}
where the sum is over all $G \in \mathcal{F}_{k+1}(P)$ containing $F$ in their boundary, and therefore
\begin{displaymath}
\sum_{F \in \mathcal{F}_k(P)} \langle  [v_F],\omega\rangle A_{n-1}(-\partial n_F,J) =
\sum_{G \in \mathcal{F}_{k+1}(P)} A_{n-1}(-n_{G},J) \sum_{F \subset G}\mathrm{sign}(F,G)\langle
 [v_F],\omega\rangle,
\end{displaymath}
but the internal sum is obviously zero by the case $k=n-1$. 
\endproof

\begin{Lemma} \label{lemma_symmetry}
Let $I,J \subset S^{n-1}$ be oriented, geodesically convex polytopes with $\dim I=k$, $\dim J=n-1-k$, such that $J \cap
E_I=\emptyset$. Then
\begin{displaymath}
A_{n-1}(\partial I,J)= (-1)^k A_{n-1}(I,\partial J).
\end{displaymath}
\end{Lemma}

\proof
Let us consider $I,J$ as singular cycles, such that the singular boundary operator on an oriented point equals its sign. Denote $\partial I=\sum_i I_i$,
$\partial J=\sum_j J_j$,  where $I_i,J_j$ are geodesically convex. Choose any point $x \in S^{n-1}$ outside $I
\cup J$, and let $H=S^{n-1} \setminus\{x\}$. Choose a form $\beta \in \Omega^{n-2}(H)$
such that $d\beta=\vol_{n-1}$. Then 
\begin{displaymath}
\partial \conv(I_{i},J)=\conv(\partial I_{i},J)+ (-1)^k\sum_j  \conv(I_{i},J_{j}),
\end{displaymath}
and since $\partial^{2}=0$, we can write 
\begin{displaymath}
\partial\sum_i \conv(I_{i},J)= (-1)^k \sum_{i,j} \conv(I_{i},J_{j}).
\end{displaymath}
Similarly, 
\begin{displaymath}
\partial \sum_j \conv(I,J_{j})=\sum_{i,j} \conv(I_{i},J_{j}).
\end{displaymath}
Therefore
\begin{align*}
A_{n-1}(\partial I,J) & = \left\langle \vol_{n-1}, \sum_i \conv(I_{i},J)\right\rangle\\
&  = \left\langle \beta,\partial \sum_i \conv(I_{i},J)\right\rangle\\
& =  (-1)^k\left\langle\beta,\sum_i \sum_j \conv(I_{i},J_{j})\right\rangle\\
& = (-1)^k A_{n-1}(I,\partial J),
\end{align*}
concluding the proof.
\endproof

The following proposition is the main result of this section. It shows that in fact $A_{n-1}$ can be uniquely extended as a bilinear functional on all chains.

\begin{Proposition}  \label{prop_volume_current}
 Given $I \in \mathcal{P}(S^{n-1})$ of dimension $k$, where $0 \leq k \leq n-2$, there exists a unique $L^1$-integrable form $\omega_I \in \Omega^{n-k-2}(S^{n-1} \setminus I)$, such that
for any $J \in \mathcal{P}(S^{n-1})$ of dimension $(n-k+2)$ with $J \cap I=\emptyset$ one has 
\begin{equation} \label{eq_omegaI}
\int_J \omega_I=A_{n-1}(-I,J). 
\end{equation}
The current $T_I \in \mathcal{D}_{k+1}(S^{n-1})$ with 
\begin{displaymath}
\langle T_I,\phi\rangle:=\int_{S^{n-1}} \omega_I \wedge \phi, \quad \phi \in \Omega^{k+1}(S^{n-1}) 
\end{displaymath}
has the following properties:
\begin{enumerate}
 \item $T_I$ is additive (i.e. $T_I=T_{I_1}+T_{I_2}$ whenever
$I=I_1 \cup I_2$ with geodesically convex polytopes $I_1,I_2$ such that $I_1 \cap I_2$ is a common face of $I_1$ and $I_2$);
\item the singular support of $T_I$ equals $I$;
\item If $k>0$, then 
\begin{equation} \label{eq_boundary_T_I}
\partial T_I= (-1)^{nk+1} \vol(S^{n-1})[[I]]+    (-1)^{n+1} T_{\partial
I},
\end{equation}
where $[[I]] \in \mathcal{D}_k(S^{n-1})$ is the $k$-current of integration over $I$;
\item If $k=0$, 
\begin{displaymath}
 \partial T_I=-\vol(S^{n-1})[[I]]+(-1)^nT_{\partial I}, 
\end{displaymath}
where we adopt the convention $\omega_{\pm 1}:=\mp \vol_{n-1} \in \Omega^{n-1}(S^{n-1}), T_{\pm 1}=\mp [[S^{n-1}]]
\llcorner \vol_{n-1} \in \mathcal{D}_0(S^{n-1})$. Then
\eqref{eq_omegaI} holds
also for $I=\pm 1$. 
\end{enumerate}
\end{Proposition}

\proof 

\emph{Step 1.} Define for $v\in S^{n-1}$ the hemisphere $H_v:=\{p\in S^{n-1}: \langle p,v\rangle>0\}$. 

Let $W \subset (S^{n-1})^n$ be the set of $n$-tuples $(p_1,...,p_n)$ belonging to some $H_v$. Define $F: W \to \R$ by $F(p_1,...p_n)=\vol_{n-1}(\Delta(p_1,...,p_n))$, the oriented volume of the geodesic simplex $\Delta(p_1,\ldots,p_n)$ with vertices $p_1,\ldots,p_n$. $F$ is well defined and smooth, since all $p_j$ lie in one hemisphere.

For two non-antipodal points $q,p \in S^{n-1}$, we define 
\begin{displaymath}
\omega_{q,p} \in \Lambda^k T ^*_qS^{n-1} \otimes
\Lambda^{n-k-2} T ^*_qS^{n-1} 
\end{displaymath}
by setting, for $u_1,\ldots,u_k \in T_q S^{n-1}$, $v_1,\ldots,v_{n-k-2}\in T_p S^{n-1}$

\begin{multline*}
\omega_{q,p}(u_1,\ldots,u_k,v_1,\ldots,v_{n-k-2})\\
:=\left.\frac{d^{n-2}}{ds^k dt^{n-k-2}}\right|_{s,t=0}
F(q, \gamma_1(s)
,\ldots,\gamma_k(s),p,\delta_1(t),\ldots,\delta_{n-k-2}(t)),
\end{multline*}
where $\gamma_i$ resp. $\delta_j$ are any smooth curves through $q$ resp. $p$ such that $\gamma_i'(0)=u_i$,
$\delta_j'(0)=v_j$. It is immediate that the definition is independent of the choice of such curves,
and that $\omega_{q,p}$ defines a unique element $\omega\in \Omega^{k,n-k-2}(S^{n-1}\times S^{n-1}\setminus \overline\Delta)$, where $ \overline\Delta=\{(q,-q):q\in S^{n-1}\} \subset S^{n-1} \times S^{n-1}$ denotes the skew-diagonal.
Note that 
\begin{multline*}
F(q,\gamma_1(\epsilon),\ldots,\gamma_k(\epsilon),p,\delta_1(\epsilon),\ldots,\delta_{n-k-2}(\epsilon))\\
=\frac{1}{k!(n-k-2)!}\omega_{q,p}(u_1,\ldots,u_k,v_1,\ldots,v_{n-k-2})\epsilon^{n-2}+o(\epsilon^{n-2}). 
\end{multline*}

Given an oriented geodesic $k$-dimensional polytope $I\subset S^{n-1}$, define $\omega_I \in \Omega^{n-k-2}(S^{n-1}\setminus I)$ by
\begin{displaymath}
\omega_I\big|_p=\int_I\omega_{-q,p}dq. 
\end{displaymath}

Let us verify that $\int_J \omega_I=A_{n-1}(-I,J)$ for an $(n-k-2)$-dimensional geodesic polytope $J$ such that $J\cap I=\emptyset$. Since both sides are additive in both $I,J$, we may assume that $I,J$ are geodesic simplices. We may further assume that there are vector fields $U_1,\ldots,U_k$ on $-I$ that are orthonormal and tangent to $-I$, and $V_1,\ldots,V_{n-k-2}$ on $J$ orthonormal and tangent to $J$. These vector fields define flow curves on $-I,J$. 

For $\epsilon>0$ one can use those curves to define a grid on $-I$, resp. $J$ denoted $\{-q_i\}$ resp. $\{p_j\}$, defining parallelograms $-Q_i$ resp. $P_j$ of volumes $\epsilon^k+o(\epsilon^k)$ resp. $\epsilon^{n-k-2}+o(\epsilon^{n-k-2})$. Note that the volume of the convex hull of two $\epsilon$-simplices is equal, up to $o(\epsilon^{n-2})$, to  $\frac{1}{k!(n-k-2)!}$ times the volume of the convex hull of the corresponding parallelograms. Thus the total volume is given by 
\begin{align*}
 A(-I,J) & =\sum A(-Q_i,P_j)\\
& =\sum _{i,j}\left(\omega_{-q_i,p_j}(U_1,\ldots,U_k,V_1,\ldots,V_{n-k-2})\epsilon^{n-2}+o(\epsilon^{n-2})\right)\\
& =\int_{(-I)\times J}\omega +o(1).
\end{align*}
Taking $\epsilon\to 0$, this proves the claim.

\emph{Step 2. } 
Let $J \subset S^{n-1}\setminus I$ be a geodesic polytope
 of dimension $(n-k-1)$. Then 
\begin{align*}
\int_J d\omega_I & = \int_{\partial J} \omega_I\\
& = A_{n-1}(-I,\partial J) \\
& =  (-1)^k A_{n-1}(-\partial I,J) \quad \text{ by Lemma \ref{lemma_symmetry}} \\
& = (-1)^k \int_J \omega_{\partial I}.
\end{align*}
It follows that on $S^{n-1} \setminus I$ we have 
\begin{displaymath}
d\omega_I= (-1)^k \omega_{\partial I}.  
\end{displaymath}

\emph{Step 3.} Let us verify that $\omega_I$ is an integrable section of $\Omega^{n-k-2}(S^{n-1} \setminus I)$, and
therefore admits a unique extension to all of $S^{n-1}$ as a current of finite mass.

Introduce spherical coordinates
\begin{align*}
 \Phi_n:[0,2\pi] \times [0,\pi]^{n-2} & \to S^{n-1}\\
(\theta_0,\theta_1,\ldots,\theta_{n-2}) & \mapsto \Phi_n(\theta_0,\theta_1,\ldots,\theta_{n-2}),
\end{align*}
which are inductively defined by 
\begin{align*}
 \Phi_2(\theta_0) & :=(\cos \theta_0,\sin \theta_0),\\
 \Phi_n(\theta_0,\theta_1,\ldots,\theta_{n-2}) & := \left(\sin \theta_{n-2} \Phi_{n-1}(\theta_0,\theta_1,\ldots,\theta_{n-3}),\cos \theta_{n-2}\right).
\end{align*}

Note that $\theta_{n-2}$ is defined on the whole sphere $S^{n-1}$ and smooth outside $\{\theta_{n-2}=0,\pi\}=S^0$, while for $i>0$, 
$\theta_{n-2-i}$ is undefined in $\{\theta_{n-1-i}=0,\pi\} \cup \{\theta_{n-1-i}\text{ undefined}\}=S^{i}$, and constitutes a coordinate outside $\{\theta_{n-2-i}=0,\pi\}$.

The volume form of $S^{n-1}$ is given by  
\begin{displaymath}
\vol_{n-1}=\prod_{i=0}^{n-3} \sin^{n-2-i} \theta_{n-2-i} \bigwedge_{i=0}^{n-2} d\theta_{n-2-i}.
\end{displaymath}

Define for $0 \leq i\leq n-2$ the vector fields 
\begin{displaymath}
X_{n-2-i}=\frac{1}{\prod_{j=0}^{i-1}\sin^{n-2-j}\theta_{n-2-j}}\frac{\partial}{\partial\theta_{n-2-i}}.
\end{displaymath}
The vector field $X_{n-2-i}$ is well defined outside the set $\{\theta_{n-2-i}=0,\pi\}$. Whenever two such vector fields
are defined, they are pairwise orthonormal. Now $\omega_I$ is integrable if
\begin{displaymath}
 \int_{S^{n-1}} |\omega_I(X_{i_1},\ldots,X_{i_{n-k-2}})| \vol_{n-1}<\infty
\end{displaymath}
for all $i_1,\ldots,i_{n-k-2}$. Let $j_1,\ldots,j_{k+1}$ be the indices not appearing in $\{i_1,\ldots,i_{n-k-2}\}$. We
consider the common level sets $C=C(\theta_{j_1},\ldots,\theta_{j_{k+1}})$, with volume element
$\sigma_C$, so that 
\begin{displaymath}
\vol_{n-1}=\left(\prod_{l=1}^{k+1}\sin^{j_l}\theta_{j_l}\wedge_{l=1}^{k+1} d\theta_{j_l}\right)\wedge \sigma_C.
\end{displaymath}
Then 
\begin{align*}
\int_{S^{n-1}} & |\omega_{I}(X_{i_1},\ldots,X_{i_{n-k-2}})| \vol_{n-1} \\
&
=\int_{\theta_{j_1},\ldots,\theta_{j_{k+1}}} \prod_{l=1}^{k+1} \sin^{j_l}\theta_{j_l} \prod_{l=1}^{k+1}d\theta_{j_l}
\int_{C(\theta_{j_1},\ldots,\theta_{j_{k+1}})}|\omega_I(X_{i_1},\ldots,X_{i_{n-k-2}})| \sigma_C.
\end{align*}

While $C(\theta_{j_{1}},...,\theta_{j_{k+1}})$ is not a geodesic polytope in $S^{n-1}$, it nevertheless holds by the
definition of $\omega_I$ that the internal integral is bounded by the total area of the sphere. Thus the entire integral
is finite. We can therefore define the current $T_I \in \mathcal{D}_{k+1}(S^{n-1})$ by 
\begin{displaymath}
 \langle T_I,\phi\rangle :=  \int_{S^{n-1} \setminus I} \omega_I \wedge \phi, \quad \phi \in \Omega^{n-k-2}(S^{n-1}).
\end{displaymath}

\emph{Step 4.  }

We prove \eqref{eq_boundary_T_I} by induction on $k$. For the induction base $k=0$, recall that $T_1=-[[S^{n-1}]]
\llcorner \vol_{n-1},
\omega_1=-\vol_{n-1}$. 

The $0$-dimensional geodesic polytope $I$ is just a point, which we may suppose to be positively oriented.
Then $-I$ is the positively oriented antipodal point. Let $S^{n-1}_\epsilon$ be the sphere $S^{n-1}$ minus the
geodesic ball of radius $\epsilon$ centered at $I$. For $g \in
C^{\infty}(S^{n-1})$, we compute
\begin{align*}
\langle \partial T_I,g\rangle & = \langle T_I,dg\rangle\\
& = \int_{S^{n-1} \setminus \{I\}}  \omega_I \wedge dg\\
& = \lim_{\epsilon \to 0} \int_{S^{n-1}_\epsilon} \omega_I \wedge dg\\
& = \lim_{\epsilon \to 0} \left[ (-1)^{n+1}\int_{S^{n-1}_\epsilon}  d\omega_I \wedge g+ (-1)^n\int_{\partial S^{n-1}_\epsilon}  \omega_I \wedge g\right]\\
& =  (-1)^{n+1} \int_{S^{n-1}} g \vol_{n-1}+ (-1)^n \lim_{\epsilon \to 0} \int_{\partial S^{n-1}_\epsilon} \omega_I \wedge g.
\end{align*}

The boundary of $S^{n-1}_\epsilon$ is an $(n-2)$-dimensional geodesic sphere around $I$. Since $\int_{\partial
S^{n-1}_\epsilon} \omega_I=A_{n-1}(-I, \partial S_\epsilon^{n-1})=(-1)^{n-1}\vol S^{n-1}_\epsilon$  (note that Lemma
(\ref{lemma_symmetry}) does not apply here, as $S_\epsilon$ is not geodesically convex), the second integral tends to
$(-1)^{n-1}g(I)$ times the volume of $S^{n-1}$.

It follows that 
\begin{displaymath}
\partial T_I= -\vol(S^{n-1}) [[I]]+  (-1)^{n} T_1,
\end{displaymath}
as claimed. 

\emph{Step 5.} 
Suppose now that $k>0$ and that \eqref{eq_boundary_T_I} holds for all polytopes of dimension strictly smaller than $k$. 
 
Define the current $U_I:= \partial T_I+(-1)^n T_{\partial I} \in
\mathcal{D}_k(S^{n-1})$. By step 2 and equation \eqref{eq_boundary_vs_differential}, $U_I$ is supported on $I$.
We have to show that $U_I=   (-1)^{nk+1} \vol(S^{n-1}) [[I]]$. 

Choose a family of closed neighborhoods $I_\epsilon$ with smooth boundary  such that $I_\epsilon$ converges to
$I$ as $\epsilon \to 0$. Define
currents $T_{I,\epsilon}, V_{I,\epsilon}, U_{I,\epsilon}$ on $S^{n-1}$ by 
\begin{align*}
 \langle T_{I,\epsilon},\phi\rangle & := \int_{S^{n-1} \setminus I_\epsilon}  \omega_I \wedge \phi, \quad \phi \in \Omega^{k+1}(S^{n-1}),\\
 \langle V_{I,\epsilon},\phi\rangle & := \int_{S^{n-1} \setminus I_\epsilon}  \omega_{\partial I} \wedge \phi,\quad \phi \in
\Omega^k(S^{n-1}),\\
 \langle U_{I,\epsilon},\phi\rangle & :=  (-1)^{n+k} \int_{\partial (S^{n-1} \setminus I_\epsilon)}  \omega_I
\wedge \phi, \quad \phi \in \Omega^k(S^{n-1}).
\end{align*}

By Step 3, $\mathbf{M}(T_{I,\epsilon}-T_I) \to 0, \mathbf{M}(V_{I,\epsilon}-T_{\partial I})\to 0$ as $\epsilon \to 0$. Since $U_{I,\epsilon}$ is given by
integration of a smooth form on a compact smooth manifold, it is a normal current (i.e. its mass and the mass of its boundary are
finite).

By Stokes' theorem and Step 2, we have  $U_{I,\epsilon}= \partial
T_{I,\epsilon} + (-1)^n V_{I,\epsilon} $. Therefore
\begin{displaymath}
 \mathbf{F}(U_{I,\epsilon}-U_I)=\mathbf{F}(\partial(T_{I,\epsilon}-T_I)+
(-1)^n(V_{I,\epsilon}-T_{\partial I})) \leq
\mathbf{M}(T_{I,\epsilon}-T_I)+\mathbf{M}(V_{I,\epsilon}-T_{\partial I}) \to 0.
\end{displaymath}

It follows that $U_I$ is a  real flat $k$-chain supported on the $k$-dimensional spherical polytope $I$.
By induction,
we have $\partial T_{\partial I}=  (-1)^{nk+n+1}\vol(S^{n-1})
[[\partial I]]$ and hence $\partial U_I=  (-1)^n \partial T_{\partial I}=  (-1)^{nk+1} \vol(S^{n-1})
[[\partial I]]$. The constancy
theorem \cite[4.1.31]{federer_book},
\cite[Proposition 4.9]{morgan_book} implies that
$U_I=  (-1)^{nk+1} \vol(S^{n-1})
[[I]]$, as claimed.

\endproof

For further use, we give the following current-theoretic interpretation of \eqref{eq_omegaI}. If $\phi \in \Omega^{k+1}_c(S^{n-1} \setminus I)$, then $T_I \cap \left([[S^{n-1}]] \llcorner \phi\right) = [[S^{n-1}]] \llcorner (\omega_I \wedge \phi)$ and hence 
\begin{displaymath}
 \langle T_I \cap  ([[S^{n-1}]] \llcorner \phi),1\rangle=\int_{S^{n-1}} \omega_I \wedge \phi=(-1)^{(n-k-2)(k+1)} \langle [[S^{n-1}]] \llcorner \phi,\omega_I\rangle.
\end{displaymath}
Let $\phi_i \in \Omega^{k+1}_c(S^{n-1} \setminus I)$ be a sequence with $[[S^{n-1}]] \llcorner \phi_i \to [[J]]$
 in $\mathcal{D}_{n-k-2,\WF([[J]])}(S^{n-1})$ Since $\WF(T_I) \cap s \WF([[J]])=\emptyset$ by the
second item, $\langle
T_I \cap ([[S^{n-1}]] \llcorner \phi_i),1\rangle \to \langle T_I \cap [[J]],1\rangle$, while $\langle [[S^{n-1}]]
\llcorner \phi_i, \omega_I\rangle \to \int_J \omega_I=A_{n-1}(-I,J)$. Therefore
\begin{equation} \label{eq_intersection_volume_current}
  \langle T_I \cap [[J]],1\rangle = (-1)^{n(k+1)} A_{n-1}(-I,J).
\end{equation}

\begin{Proposition} \label{prop_boundary_tildet}
 Let $P$ be a polytope and let $(T,C):=E(M([P]))$ the associated currents. Decompose $T=t+T'$ as in Proposition
\ref{prop_def_convolution}, where $t$ is the $(0,n-1)$-component of $T$ and $T'=\sum_{k=1}^{n-1} \sum_{F \in
\mathcal{F}_k(P)} A_{v_F,\check{n}(F,P)}$. Let 
\begin{displaymath}
 \tilde T:= \frac{1}{\vol(S^{n-1})} \sum_{k=1}^{n-1}  (-1)^{nk+k+1}\sum_{F \in \mathcal{F}_k(P)}   [v_F]  \times
T_{\check n(F,P)}
\in \mathcal{D}_n( SV).
\end{displaymath}
Then 
\begin{displaymath}
 \partial \tilde T=T'. 
\end{displaymath}
\end{Proposition}

\proof  
By \eqref{eq_boundary_product} and Proposition \ref{prop_volume_current} we have
\begin{align*}
\partial \tilde{T} & = \frac{1}{\vol(S^{n-1})}  \sum_{k=1}^{n-1}  (-1)^{nk+k+1} \sum_{F \in \mathcal{F}_k(P)} 
\partial([v_F] \times T_{\tilde
n(F,P)})\\
& = \frac{1}{\vol(S^{n-1})}  \sum_{k=1}^{n-1}  (-1)^{nk+1}  \sum_{F\in \mathcal{F}_k(P)} [v_F]
\times \partial
T_{\tilde
n(F,P)}\\
& =\sum_{k=1}^{n-1} \sum_{F \in \mathcal{F}_k(P)} A_{v_F,\check{n}(F,P)}+\frac{1}{\vol(S^{n-1})} \sum_{k=1}^{n-1} \sum_{F\in
\mathcal{F}_k(P)} (-1)^{nk+n} [v_F]  \times T_{\partial n_F}.
\end{align*}

Let $ 1 \leq k \leq n-1$ be fixed and let $J \subset S^{n-1}$ be a $(k-1)$-dimensional geodesic polytope not
intersecting any $\partial n_F$ for $F \in \mathcal{F}_k(P)$. Lemma \ref{lemma_cancellation} implies that  
\begin{displaymath}
\sum_{F \in \mathcal{F}_k(P)} \langle [v_F],\omega\rangle T_{\partial n_F} \cap [[J]]
=0
\end{displaymath}
for all $\omega \in \Omega_c^k(\R^n)$. It follows that $\sum_{F \in \mathcal{F}_k(P)}  [v_F]  \times T_{\partial
n_F}=0$, and the statement follows.
\endproof

\section{Compatibility of the algebra structures}
\label{sec_compatibility}

From now on, we fix an orientation and a Euclidean scalar product on $V$ and identify $\Dens(V) \cong
\mathbb{C}, \ori(V) \cong \mathbb{C}, S^*V \cong SV=V \times S^{n-1}$. It then holds that 
\[\langle \ast_1 T,\omega\rangle =(-1)^{nk+nl+k}\langle T,\ast_1\omega \rangle \]  for $T\in \mathcal D_{k,l}(S^*V)^{tr}$, $\omega \in \Omega^{k,n-1-l}(S^*V)^{tr}$.

\begin{Proposition} \label{prop_star_product}
Let $v_i \in \Lambda^{k_i}_sV$ and $T_i \in \mathcal{D}_{l_i}(S^{n-1}), i=1,2$ be currents on the sphere such that
$T_1 \cap T_2 \in
\mathcal{D}_{l_1+l_2-n+1}(S^{n-1})$ is defined. Then 
\begin{displaymath}
 \ast_1 ([v_1] \times T_1) \cap \ast_1([v_2] \times T_2) =  (-1)^{k_1(n-k_2-l_2-1)}\ast_1 \left([v_1
\wedge v_2] \times (T_1
\cap T_2)\right).
\end{displaymath}
In particular, for $A_{v_i,N_i} \in Y_{k_i}$ with $N_1$ and $N_2$ being transversal, one has
\begin{displaymath}
\ast_1 A_{v_1,N_1} \cap \ast_1 A_{v_2,N_2}= \ast_1 A_{v_1 \wedge v_2,N_1 \cap N_2}.
\end{displaymath}
\end{Proposition}

\proof
Let $v \in \Lambda^k_s V, T \in \mathcal{D}_l(S^{n-1})$. We compute  for $\gamma_1 \in \Omega_c^{n-k}(V),
\gamma_2 \in \Omega_c^l(S^{n-1})$ 
\begin{align*}
\langle \ast_1 ([v] \times T),\pi_1^* \gamma_1 \wedge \pi_2^*\gamma_2\rangle
& =(-1)^{nk+nl+k}\langle [v] \times T,\ast_1 (\pi_1^* \gamma_1 \wedge \pi_2^*\gamma_2)\rangle\\
& =(-1)^{nk+nl+k+\binom{k}{2}} \langle [v],\ast \gamma_1\rangle \cdot \langle T,\gamma_2\rangle\\
& =(-1)^{nk+nl+k+\binom{k}{2}+k(n-k)} \langle [\ast v], \gamma_1\rangle \cdot \langle T,\gamma_2\rangle\\
& =(-1)^{\binom{k}{2}+nl} \langle [\ast v] \times T,\pi_1^* \gamma_1 \wedge \pi_2^*\gamma_2\rangle,
\end{align*}
i.e.
\begin{equation} \label{eq_star1}
 \ast_1 ([v] \times T)=(-1)^{\binom{k}{2}+nl} [\ast v] \times T.
\end{equation}

Let now $v_i \in \Lambda^{k_i}_s V, T_i \in \mathcal{D}_{l_i}(S^{n-1}), i=1,2$. We have $[*v_1] \cap [*v_2]=[*(v_1 \wedge v_2)]$. Using \eqref{eq_intersection_product} it follows that  
\begin{displaymath}
 ([*v_1] \times T_1) \cap ([*v_2] \times T_2) =  (-1)^{k_1(n-1-l_2)} [*(v_1 \wedge v_2)]
\times (T_1 \cap T_2).
\end{displaymath}
The statement now follows from \eqref{eq_star1}.
\endproof

\begin{Definition}
Two elements $x,y \in \Pi(V)$ are in general position,
if any two normal cones to a pair of faces of $x$ and $y$ are transversal.
\end{Definition}

\begin{Proposition}
Given two elements $x,y\in\Pi(V)$ in general position, the convolution
$M(x) \ast M(y)$ is well-defined, and $M(x\cdot y)=M(x) \ast M(y)$.
\end{Proposition}

\proof
By linearity, it suffices to consider $x=[P]$ and $y=[Q]$ for some
polytopes $P,Q$. 

The wavefront of the current $A_{v,N} \in Y_k$ is the conormal bundle to $N$ in $S^{n-1}$. Thus $\WF(
A_{v_1,N_1}) \cap
\WF( A_{v_2,N_2})=\emptyset$ if and only if for all $x \in N_1 \cap N_2$ we have 
\begin{displaymath}
T_x N_1+T_x N_2=T_x S^{n-1},
\end{displaymath}
that is, if and only if $N_1$ and $N_2$ are transversal.

Recall from Lemma \ref{lemma_image_of_polytope} that for $0 \leq k \leq n-1$ 
\begin{displaymath}
T(M_{n-k}[P])=\sum_{F \in \mathcal{F}_k(P)} A_{v_F,\check{n}(F,P)}.
\end{displaymath}
Thus, given $[P],[Q] \in \Pi(V)$ in general position, the normal cones
in $S^{n-1}$ have disjoint wavefronts. 

Let $F$ be a face of $P$ and $G$ a face of $Q$. If $\dim F+\dim G \geq n$ then $\check{n}(F,P) \cap
\check{n}(G,Q)=\emptyset$
by transversality. 

Thus, by Proposition \ref{prop_star_product}
\begin{align*}
\ast_1 T(M([P]) \ast M([Q])) & =\sum_{\dim F+\dim G<n} \ast_1 A_{v_F,\check{n}(F,P)} \cap \ast_1
A_{v_G,\check{n}(G,Q)}\\
& =\sum_{\dim F+\dim G<n} \ast_1 A_{v_F \wedge v_G,\check{n}(F,P) \cap \check{n}(G,Q)}\\
& =\ast_1 \sum_{H \in
\mathcal{F}(P+Q)} A_{v_H,\check{n}(H,P+Q)},
\end{align*}
and so 
\begin{displaymath}
T(M[P] * M[Q]) =\sum_{H\in \mathcal{F}(P+Q)} A_{v_{H},\check{n}(H,P+Q)}=T(M[P+Q]).
\end{displaymath}

It remains to verify that $C(M[P] \ast M[Q])=C(M[P+Q])$.  We set 
\begin{displaymath}
T_1=t_1+T_1':=T(M[P]), T_2=t_2+T_2':=T(M[Q]) 
\end{displaymath}
as in Proposition \ref{prop_def_convolution}. Then $t_1=\sum_{F \in \mathcal{F}_0(P)} A_{F,n(F,P)}=\pi^*( [[V]]
\llcorner
\vol_n)$, hence $\alpha_1=1$. Similarly $\alpha_2=1$. 

Let $F$ be a $k$-dimensional face of $P$ and $G$ an $(n-k)$-dimensional face of $Q$,  where $1 \leq k \leq n-1$.  By Proposition
\ref{prop_star_product},
\begin{align*}
\ast_1([v_F] \times T_{\check n(F,P)})\cap \ast _1 (A_{G, \check n(G, Q)}) & =
\ast_1([v_F\wedge v_G]
 \times (T_{\check n(F,P)}\cap [[\check n(G,Q)]]).
\end{align*}

Summing over all faces of $P$ and $Q$, using \eqref{eq_intersection_volume_current} and Proposition \ref{prop_boundary_tildet}, we obtain that 
\begin{align*}
C( & M[P] *  M[Q]) =  \pi_* *_1^{-1} (*_1 \tilde T_1 \cap *_1 T_2') + \alpha_1 C_2+\alpha_2 C_1\\
 & =   \sum_{k=1}^{n-1} (-1)^{nk+k+1} \sum_{F \in \mathcal{F}_k(P)}  \pi_* *_1^{-1} (*_1 ([v_F] \times
T_{\check n(F,P)}) \cap *_1 T_2') \\
& \quad + \alpha_1 C_2+\alpha_2 C_1\\
& =\sum_{k=1}^{n-1}  (-1)^{nk+k+1} \sum_{\substack{F \in \mathcal{F}_k(P)\\G \in \mathcal{F}_{n-k}(Q)}}  \pi_*
\left([v_F\wedge v_G]
\times (T_{\check n(F,P)}\cap [[\check n(G,Q)]])\right)  \\
& \quad + \vol(P)[[V]]+\vol(Q)[[V]]\\
& =\sum_{k=1}^{n-1}   (-1)^{n(n-k)+nk+k+1} \sum_{\substack{F \in \mathcal{F}_k(P)\\G \in
\mathcal{F}_{n-k}(Q)}} [v_F\wedge v_G] \vol(\conv(-\check n(F,P),\check n(G,Q))) \\
& \quad + \vol(P)[[V]]+\vol(Q)[[V]]\\
& =\sum_{k=1}^{n-1}  (-1)^{n+k+1} \sum_{\substack{F \in \mathcal{F}_k(P)\\G \in \mathcal{F}_{n-k}(Q)}} 
[v_F\wedge v_G] \vol(\conv(-\check n(F,P),\check n(G,Q))) \\
& \quad + \vol(P)[[V]]+\vol(Q)[[V]].
\end{align*}

Note that the pair $(-1)^{n+\dim F+1}(v_F,-\check n(F,P))$ is positively oriented, and coincides with
$(v_{-F}, \check
n(-F,-P))$. It follows by 
\cite[Eq. 5.66]{schneider_book14} that 
\begin{align*}
C(M[P] * M[Q]) & = \sum_{k=1}^{n-1} \sum_{\substack{F \in \mathcal{F}_k(P)\\G \in \mathcal{F}_{n-k}(Q)}} [v_{-F}
\wedge v_G] \vol(\conv(\check n(-F,-P),\check n(G,Q))) \\
& \quad + \vol(P)[[V]]+\vol(Q)[[V]] \\
& = \vol(P +Q)[[V]]  \\
& = C(M[P+Q]). 
\end{align*}
This finishes the proof. 
\endproof

\def\cprime{$'$}

\end{document}